\documentclass[12pt,onecolumn,draftcls,journal,letterpaper]{IEEEtran}
  \usepackage[dvips]{graphicx}
   \graphicspath{{figs/}}

\usepackage{dblfloatfix}    

\usepackage[cmex10]{amsmath}

\usepackage{array}

\usepackage[caption=false,font=footnotesize]{subfig}

\usepackage{url}

\usepackage{amsmath}
\usepackage{amsfonts}
\usepackage{mathrsfs}
\usepackage{epstopdf}
\usepackage{graphicx}
\newcommand{\xb}{\textbf{x}} 
\newcommand{\ub}{\textbf{u}}

\newcommand{\subj}{\text{subj. to}}

\usepackage{tikz}
\usetikzlibrary{arrows}
\usetikzlibrary{spy}
\usetikzlibrary{calc}

\tikzset{dashdot/.style={dash pattern=on .4pt off 3pt on 4pt off 3pt}}

\usetikzlibrary{plotmarks}
\usepackage{pgfplots}

\begin{document}
%
\title{Optimal Rendezvous Trajectory for Unmanned Aerial-Ground Vehicles}

\author{A.~Rucco,
	P.B.~Sujit, 
        A.P.~Aguiar,
        J.B.~Sousa,
        and~F. L.~Pereira
\thanks{A. Rucco is with the Department of Engineering, Universit\`a del Salento, Lecce, 73100, Italy, e-mail: 
{\tt alessandro.rucco@unisalento.it}. }
\thanks{P.B. Sujit is with the Department of Electrical and Computer Engineering, Faculty of Engineering and Department of Electronics and Communications Engineering, IIIT Delhi, India, 110020, e-mail: 
{\tt sujit@iiitd.ac.in}. }
\thanks{A.P. Aguiar, J.B. Sousa and F. L. Pereira are with Department of Electrical and Computer Engineering, Faculty of Engineering, University of Porto, 4200-465, Portugal,  e-mail: 
{\tt \{pedro.aguiar,~jtasso,~flp\}@fe.up.pt}.
}
}

\maketitle

\begin{abstract}
	Fixed-wind unmanned aerial vehicles (UAVs) are essential for low cost aerial surveillance and mapping applications in remote regions. One of the main limitations of UAVs is limited fuel capacity and hence requires  periodic refueling to accomplish a mission. The usual mechanism of commanding the UAV to return to a stationary base station for refueling  can result in fuel wastage and inefficient mission operation time. Alternatively, unmanned gound vehicle (UGV) can be used as a mobile refueling unit where the UAV will rendezvous with the UGV for refueling. In order to accurately perform this task in the presence of wind disturbances, we need to determine  an optimal trajectory in 3D taking UAV and UGV dynamics and kinematics into account. In this paper, we propose an optimal control formulation to generate a tunable UAV trajectory for rendezvous on a moving UGV taking wind disturbances into account. 
	By a suitable choice of the value of an aggressiveness index in our problem setting, we are able to control the UAV rendezvous behavior. Several numerical results are presented to show the reliability and effectiveness of our approach. 
\end{abstract}

\IEEEpeerreviewmaketitle

\section{Introduction}
\begin{figure*}
	\centering
	\subfloat[]{\label{fig:scenario1}\includegraphics[width=7.5cm]{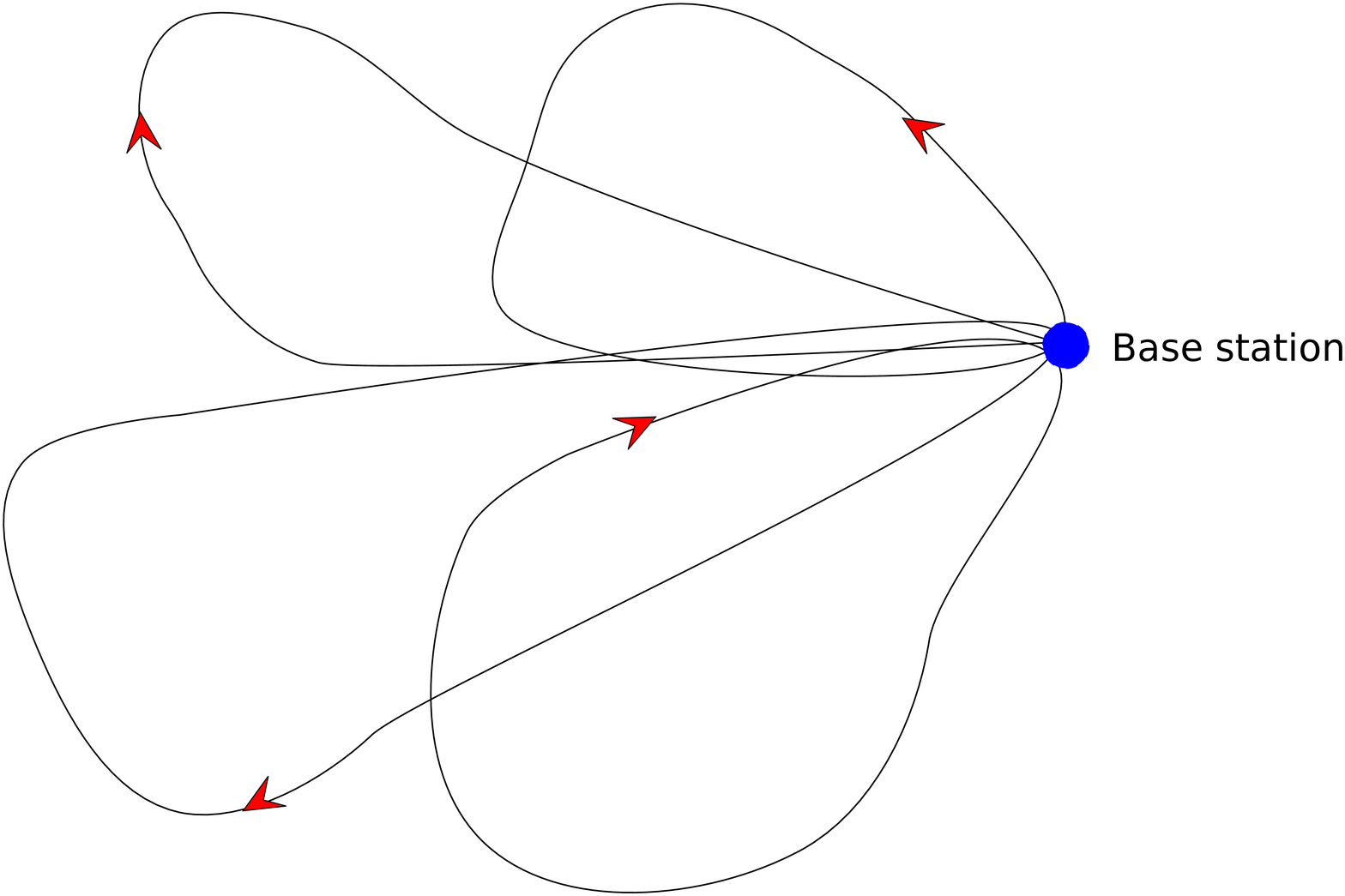}}
	\subfloat[]{\label{fig:scenario2}\includegraphics[width=7.5cm]{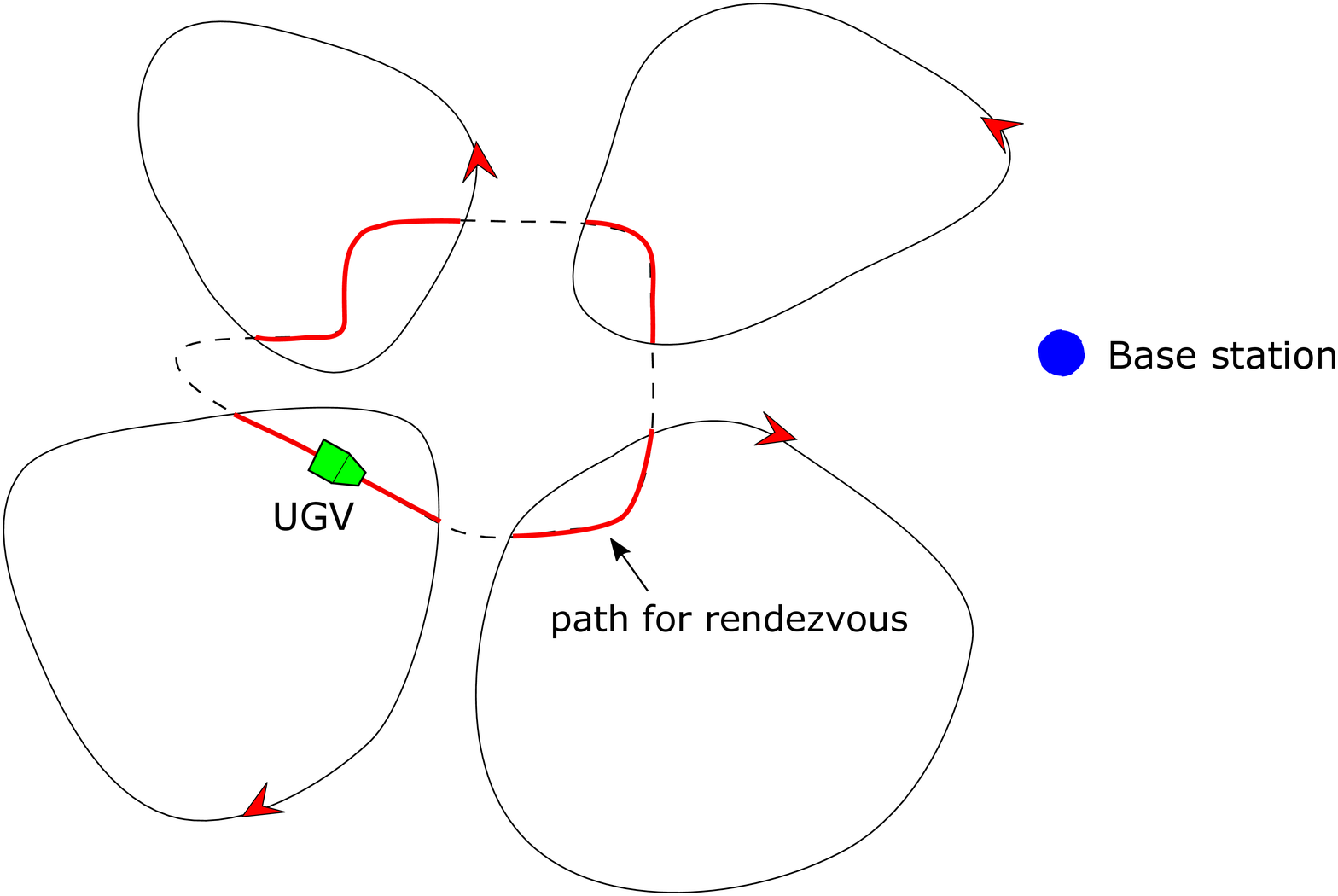}}
	\caption{(a) A field deployment where UAVs visit a base station located at a distant for refueling. (b) A UGV is deployed for refueling with a predefined UGV path and the time for rendezvous.} \label{fig:scenario}
\end{figure*}

\IEEEPARstart{F}{ixed}-wing unmanned aerial vehicles (UAVs) are essential components of remote monitoring applications like surveillance, mapping, aerial photography, etc., where the UAVs need to cover large regions. 
Typical UAVs used for these applications are of low cost with limited fuel capacity and hence require periodic refueling to accomplish the mission.
For the case of using low cost UAVs, these ones have however limited fuel capacity and require periodic refueling to accomplish the mission. 
In these scenarios, airborne docking for mid-air refueling has become recently a major research area, see e.g., \cite{fravolini2004modeling, valasek2005vision}. However, the wake effects of the tanker on the UAV makes the analysis and design of the control scheme particularly challenging (e.g., a large amount of experimental data are needed). In \cite{nichols2014aerial, sun2014optimal}, a passive towed cable system is used to retrieve the UAV, thus avoiding wake phenomena. On the other hand, a robust vision tracking method is required for the UAV to overcome some hardware limitations of the vision system (mostly when the UAV gets closer to the drogue). The most simple solution is to deploy an immobile base station in a fixed location to oversee the operation and to refuel the UAVs, as shown in Figure~\ref{fig:scenario1}. The base station may be located at a distant which diminishes the utility of the UAVs fuel per mission. Instead of an immobile unit, an unmanned ground vehicle (UGV) can be deployed that can refuel the UAVs at different locations, and hence reducing the UAV refueling time, which increases the coverage area per refuel as shown in Figure \ref{fig:scenario2}. In order to accomplish this capability, there is a need to develop techniques for UAV rendezvous with the moving UGV. 

Cooperative UAV and UGV teams have been previously used for several surveillance applications. 
For instance, the UAV can provide useful information (e.g., data from aerial images) to the UGV for path planning and target detection 
\cite{Cantelli2013,Grocholsky2006,Sauter2009,Yu2011}. 
In a different application, Tokekar et al. \cite{Tokekar2013} used an UAV to acquire points of nitrogen sampling in a field and the UGV used these points to create a path of one-in-a-set. 
In this paper, we are concerned about using the UGV as a refueling mobile station and hence the UAV needs to generate a trajectory  such that it can rendezvous with the moving UGV. 

The UAV, UGV rendezvous 
can be considered either as a docking or landing problem. Aerial rendezvous between multiple aircrafts for refueling \cite{Kampoon2010,Burns2007} and formation flight \cite{McLain2000,Harl2008,Yamasaki2010} are related but the type of vehicles taken into account are the same and secondly, the rendezvous typically is in 2D, unlike the landing, which is in 3D. Carnes et al. \cite{Carnes2015} developed an auto-takeoff and auto-landing capabilities for a low-cost UAV, which is essential for many of the envisioned applications. Nonetheless, the trajectories are not optimized, which is one of the key contributions of this paper. Kim et al. \cite{Kim2013} developed a vision based net-landing controller for a UAV. The controller is based on pure-pursuit guidance law. Daly et al. \cite{Daly2011} developed a landing controller for a quad-rotor which can hover and land on a moving vehicle. However, landing using a fixed-wing aerial vehicle instead of a quad-rotor onto a moving vehicle is much more challengeable \cite{barber2009vision}. Another relevant literature is the rendezvous/landing guidance with impact angle constraints where the impact angle is the angle about which the landing or rendezvous takes place \cite{Smith2008,Ratnoo2015}. In those works, the trajectories are not optimized. 

The contributions of this paper are as follows. First, we propose an optimization-based strategy for the generation of optimal UAV rendezvous trajectory onto a moving UGV. In order to generate realistic rendezvous trajectories, the strategy has to take dynamics and kinematics of the UAV and UGV into account.  
The coupled UAV-UGV dynamics and the constraints arising from the rendezvous maneuver make the design of the strategy complex. 
We set up the rendezvous optimal control problem in terms of a suitable error dynamics which describe the coupled dynamics. 
The error dynamics make the analysis and design of the rendezvous strategy simpler, because the key for achieving successful rendezvous is that the error coordinates are zero at the rendezvous point. 
Second, we identify an aggressiveness index in our rendezvous optimal control problem which allows us to control the UAV rendezvous behavior. The aggressiveness index is based on the performance limitations of the UAV (i.e., the constraint limits on the state, input variables), thus allowing us to compute aggressive trajectories (several dynamic constraints are active while the UAV is approaching the UGV) or very smooth ones. 
The proposed optimal solution framework for the UAV-UGV rendezvous can be seen as a framework which allows one to select (in form of tuning knob) the type of UAV trajectory. 
Finally, through numerical computations, we show the effectiveness of our approach and discuss a set of interesting features of the rendezvous trajectories. 

The rest of the paper is organized as follows. In Section~\ref{sec:pf}, we propose the optimal control formulation for the UAV-UGV rendezvous. In Section~\ref{sec:solution},  we describe the optimal control based strategy for effectively solving the rendezvous optimal control problem.  This technique is evaluated through numerical computations and illustrated in Section~\ref{sec:results}. The conclusions are given in Section~\ref{sec:conclude}.

\begin{table}[ht]   \caption{Nomenclature} 
  \centering 
\begin{tabular}{l l} 
  $p_i = [x_i,y_i,z_i]^T$ & \!\!\!UAV ($i = A$) and UGV ($i = G$) position, $m$\\   
  $v_{i}$ & \!\!\!UAV ($i = A$) and UGV ($i = G$) ground-speed, $m/s$\\  
  $\chi_{i}$ & \!\!\!UAV ($i = A$) and UGV ($i = G$) course angle, $rad$\\
  $\gamma_A, \phi_A, \psi_A $ & \!\!\!UAV flight path angle, roll angle, heading angle, $rad$ \\
  $v_{a}$ & \!\!\!UAV airspeed, $m/s$ \\
  $\gamma_{a}$ & \!\!\!UAV air-flight path angle, $rad$ \\
  $T$ & \!\!\!Thrust, $N$ \\
  $D$ &\!\!\!Drag force, $N$ \\
  $L$ & \!\!\!Lift force, $N$\\  
  $m$ & \!\!\!Mass, $kg$ \\
  $g$ & \!\!\!Gravitational acceleration, $m/s^2$\\
  $\rho$ & \!\!\!Air density, $kg/m^3$\\
  $S$ & \!\!\!Surface area of the wing, $m^2$\\
  $C_L$ & \!\!\!Lift coefficient \\
  $C_D$ & \!\!\!Drag coefficient \\
  $C_{D_0}$ & \!\!\!Drag coefficient at zero lift \\
  $k_{D/L}$ & \!\!\!Induced drag factor \\
  $n_{lf}$ & \!\!\!Load factor\\  
  $\alpha$ & \!\!\!Angle of attack, $rad$\\
  $a_{lon}, a_{lat}$ & \!\!\!UGV longitudinal and lateral acceleration, $m/s^2$\\
  $\sigma_G$ & \!\!\!UGV path curvature, $1/m$ \\
  $s_G$ & \!\!\!UGV path coordinate, $m$ \\ 
  $e = [e_x,e_y,e_z]^T$ & \!\!\!Longitudinal, lateral, and vertical error coordinates, $m$\\  
  $e_\chi, e_\gamma, e_\phi$ & \!\!\!Course angle, flight path, and roll error angles, $rad$\\  
  $e_v$ & \!\!\!Speed error, $m/s$\\  
  $v_w$ & \!\!\!Wind speed, $m/s$\\  
  $w_x,w_y,w_z$ & \!\!\!Wind velocity components in the inertial frame, $m/s$\\  
\end{tabular}
  \label{table:nomenclature} 
\end{table}

\section{Problem Formulation} \label{sec:pf}
In this section, we address the rendezvous problem of a fixed-wing UAV onto a moving UGV. 
We first introduce the equations of motion for UAV and UGV, and outline the constraints. 
Second, we describe the UAV and UGV dynamics with respect to a suitable error dynamics (i.e., the velocity frame of the UGV).  
We then formulate the rendezvous problem with respect to the coupled UAV-UGV dynamics.  
In Table~\ref{table:nomenclature} we provide a list of the symbols used in
the paper. 

\subsection{UAV dynamic model}
\begin{figure}[htpb]
  \begin{center}
    \subfloat[Wind triangle projected onto the $x-y$ plane.]
    {\includegraphics[width=6cm]{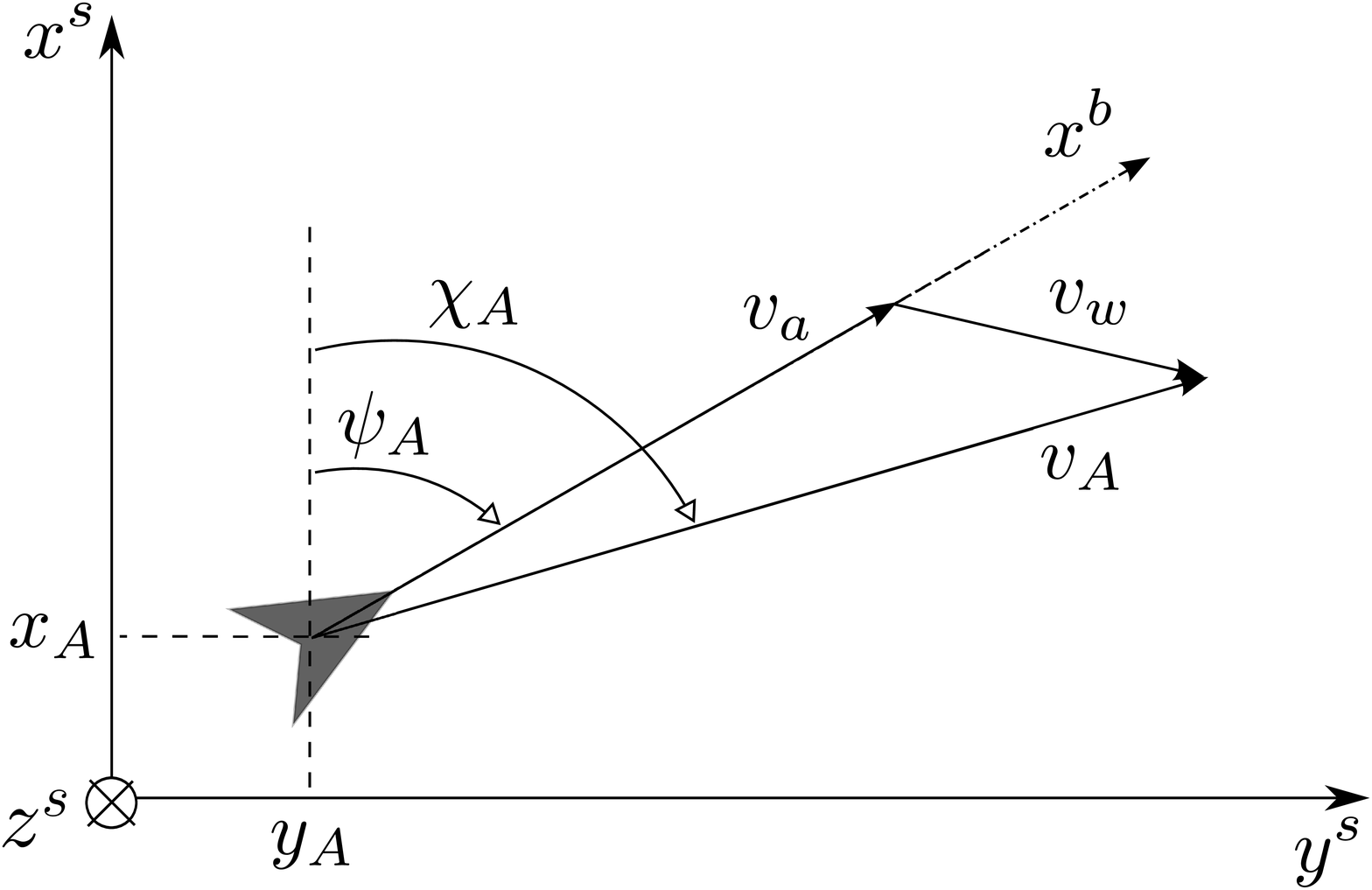} }  \,\,\,
    \subfloat[Wind triangle projected onto the $x-z$ plane.]
    {\includegraphics[width=6cm]{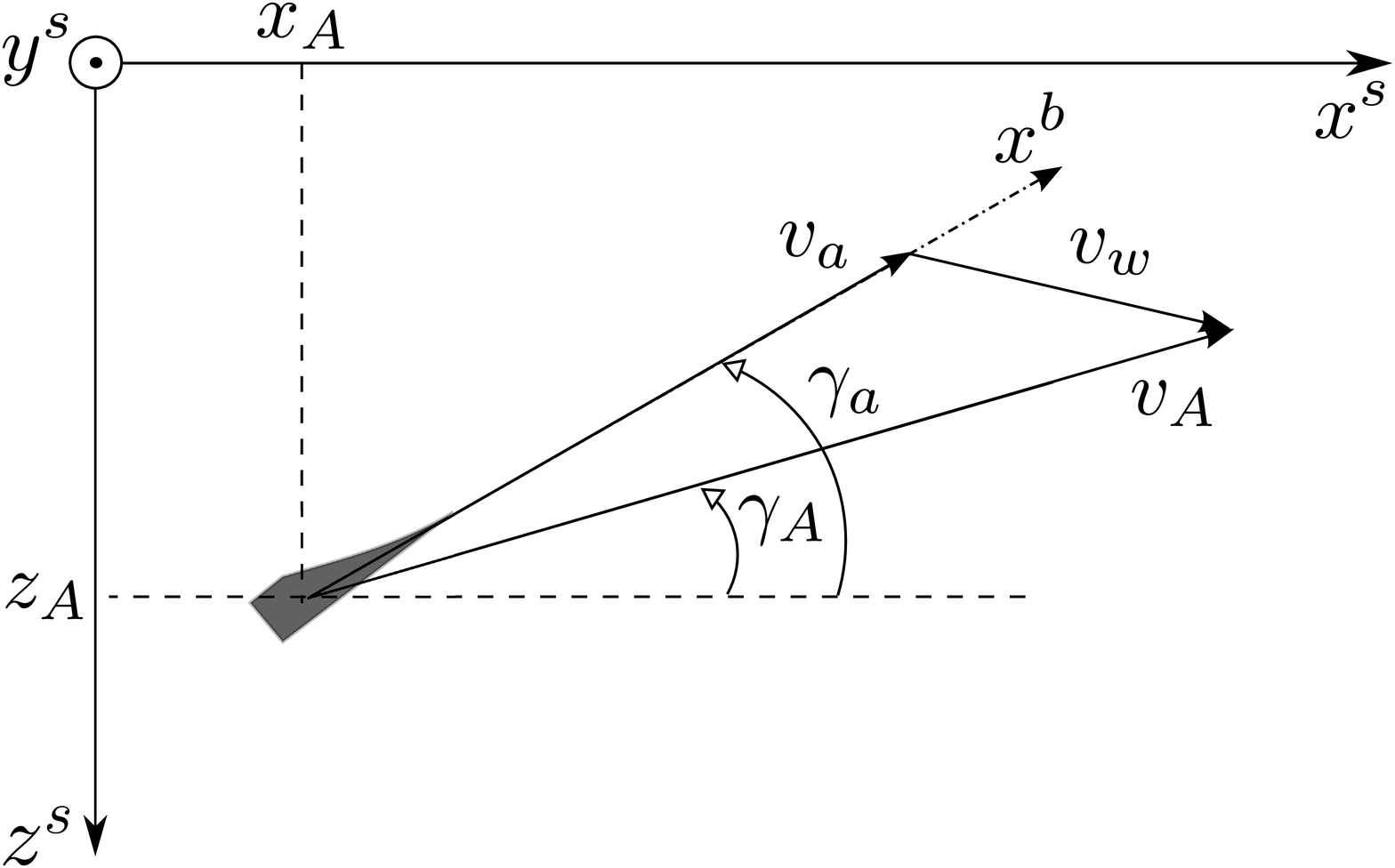} }
    \caption{The wind triangle. }
    \label{fig:WindTriangle}
  \end{center}
\end{figure}

We use a 3D point mass model for the aerial vehicle \cite{beard2012small}. The six DOF equations of motion can be written as 
where, $L = \frac{1}{2}\rho v_a^2SC_L, D = \frac{1}{2}\rho v_a^2SC_D,$ and $C_D=C_{D_0} + K_{D/L} C_L^2$. 
The airspeed, $v_a$, and the ground speed, $v_A$, are related by  
\begin{equation}  \label{eq:vgchigamma_vapsgamma}
  \begin{split}
      v_A\cos \chi_A \cos \gamma_A &= v_a\cos \psi_A \cos \gamma_a + w_x, \\
      v_A\sin \chi_A \cos \gamma_A &= v_a\sin \psi_A \cos \gamma_a + w_y, \\
      -v_A\sin\gamma_A &= -v_a \sin \gamma_a + w_z.
  \end{split} 
\end{equation}
where $w_x$, $w_y$ and $w_z$ are the wind components in the inertial frame. 
Exploiting the wind triangle, see Figure~\ref{fig:WindTriangle}, the airspeed, the heading angle, and the air-mass-referenced flight path angle are given by
\[
v_a \!=\! \sqrt{ v_A^2 \!-\! 2v_A(w_x \cos{\chi_A}\cos{\gamma_A} \!+\! w_y\sin{\chi_A}\cos{\gamma_A} \!-\! w_z \sin{\gamma_A}) \!+\! v_w^2 }\, ,
\]
\[
\gamma_a = \arcsin{\left( \frac{v_A \sin \gamma_A + w_z}{v_a} \right)}\, ,
\] 
\[
\psi_A = \chi_A - \arcsin{\left( \frac{-w_x \sin\chi_A + w_y \cos\chi_A}{v_a \cos\gamma_a} \right)} \, .
\]
We consider three control inputs for the UAV: $u_1=T$, $u_2=\dot \phi_A$, and $u_3=C_L$. In particular, we vary the thrust of the vehicle to affect the airspeed of the UAV. The $u_2$ is the roll rate by which the UAV heading angle and the flight path angle are updated. The $u_3$ is the lift coefficient which we assume to operate in the linear region and hence approximately a linear function of the angle of attack $\alpha$~\cite{beard2012small}. 
\begin{equation}  \label{sys:UAVdyn_Wind}
  \begin{split}
      \dot{x}_A &= v_A\cos \chi_A \cos \gamma_A, \\
      \dot{y}_A &= v_A\sin \chi_A \cos \gamma_A, \\
      \dot{z}_A &= -v_A\sin\gamma_A, \\
      \dot{v}_A &= \frac{u_1-D}{m}-g\sin \gamma_A, \\
      \dot{\gamma}_A &= \frac{1}{v_A}\left(\frac{L\cos\phi_A}{m}-g\cos\gamma_A\right), \\
      \dot{\chi}_A &= \frac{1}{v_A\cos \gamma_A}\left(\frac{L\sin\phi_A \cos{( \chi_A - \psi_A )}}{m}\right), \\
      \dot{\phi}_A &= u_2,
  \end{split} 
\end{equation}
The UAVs have state and input constraints. 
In particular, the airspeed, $v_a$, the load factor, $n_{lf} = \frac{L}{m g}$, and the flight path angle, $\gamma_A$, are bounded by $v_{min}$ and $v_{max}$, $n_{lf\,min}$ and $n_{lf\,max}$, $\gamma_{min}$ and $\gamma_{max}$, respectively. 
The thrust is constrained to be positive and less than the maximum value $u_{1 \, max}$. 
Moreover, the roll angle, the roll rate, and the lift coefficient are bounded in module by $\phi_{max}$, $u_{2 \, max}$, $u_{3 \, max}$, respectively. 
More specifically, the following state and input constraints
are imposed on the model:
\begin{equation} \label{sys:UAVConstraints}
\begin{split}
  v_{min} \leq v_a \leq v_{max}\,, \\ 
  n_{lf \, min} \leq n_{lf} \leq n_{lf \, max}\,, \\ 
  \gamma_{min} \leq \gamma_{A} \leq \gamma_{max}\,, \\ 
  0 \leq u_1 \leq u_{1 \, max}\,, \\ 
  | \phi_A | \leq \phi_{max}\,, \\
  | u_2 | \leq u_{2\, max}\,, \\ 
  | u_3 | \leq u_{3\, max}\,. \\ 
\end{split} 
\end{equation}
The UAV parameters, aerodynamic coefficients and the constraint parameters used in the paper are given in Appendix. 

\subsection{UGV dynamic model}
We model the UGV as a 2D point mass model~\cite{bayer2012trajectory}. In this case, the equations of motion are 
\begin{equation}  \label{sys:UGVdyn}
  \begin{split}
	\dot{x}_G &= v_G\cos \chi_G \, , \\
	\dot{y}_G &= v_G\sin \chi_G \, , \\
	\dot{v}_G &= a_{lon}\, , \\
	\dot{\chi}_G &=v_G \sigma_G \, . \\
  \end{split} 
\end{equation}
We recall that the UGV can move on a pre-determined path as the one shown in Figure~\ref{fig:scenario2}. 
Therefore, we take the control input of the UGV to be the longitudinal acceleration, $u_4 = a_{lon}$. The lateral acceleration can be written as $a_{lat} = v^2_G \sigma_G$, where $\sigma_G$ is the (fixed) path curvature~\cite{velenis2008minimum}. 
Note that we describe the UGV curvature as a function of the path coordinate (or arc length coordinate) $s_G(t) = \int_0^t \sqrt{\dot{x}_G(\tau) + \dot{y}_G(\tau)} d\tau$. 
In other words, the UGV can accelerate/decelerate along the fixed path defined by the curvature. 

Due to the tire-road force interaction, the vehicle acceleration is limited by the so called friction circle (more generally friction ellipse)~\cite{rucco2015efficient}. 
Here, we take into account a circular acceleration constraint: 
the acceleration has to be less than or equal to $a_{max}$, i.e., 
\begin{equation} \label{sys:UGVConstraints}
\begin{split}
   a^2_{lon} + a^2_{lat} \leq a^2_{max}\,. \\ 
\end{split} 
\end{equation}

\subsection{Error dynamics}
\begin{figure}[htpb]
  \begin{center}
    \includegraphics[width=8cm]{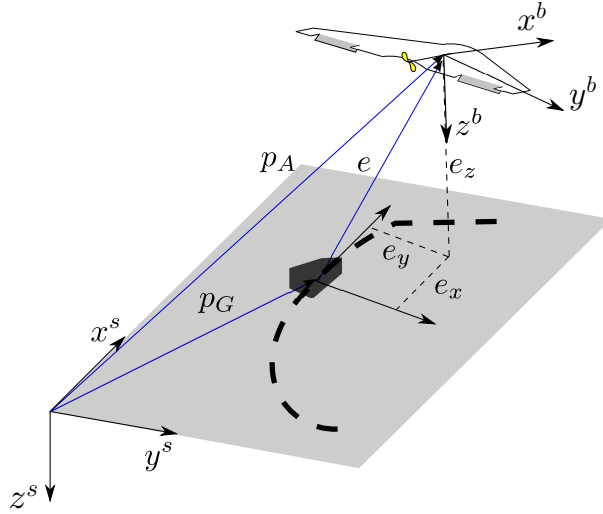}
    \caption{Error space frames and the fixed-wing UAV body frame. 
     }
    \label{fig:new_coord_3D}
  \end{center}
\end{figure}

The coordinates of the aerial vehicle expressed in the inertial frame, $p_A = [x_A, y_A, z_A]^T$, can be defined with respect to the position of the ground vehicle, $p_G = [x_G, y_G, z_G]^T$, as following, see Figure~\ref{fig:new_coord_3D}, 
\begin{equation} \label{eq:xtow} 
p_A = p_G
  +
  R_{z}({\chi}_{G}) e\,,
\end{equation}
where $e = [e_x, e_y, e_z]^T$ is the error vector expressed in the body-frame of the UGV and
\[
R_z({\chi}_{G}) = 
  \left[
    \begin{array}{ccc}
      \cos{{\chi}_{G}} & -\sin{{\chi}_{G}} & 0 \\
      \sin{{\chi}_{G}} & \cos{{\chi}_{G}} & 0\\
      0 & 0 & 1\\
    \end{array}
  \right],
\] 
is the rotation matrix transforming vectors from the error frame (i.e., the velocity frame of the UGV) into the inertial frame. 
It is worth noting that, since the altitude of UGV is constant and equal to zero ($z_G = 0$), the vertical error coordinate is equal to the altitude of the UAV, i.e., $e_z = z_A$. 

Next, we compute the expression of $\dot{e} = [\dot{e}_x, \dot{e}_y, \dot{e}_z]^T$. By  differentiating \eqref{eq:xtow} with respect to the time $t$, we get
\begin{equation} \label{eq:noname}
\dot{p}_A =  \dot{p}_G
  + \left[
    \begin{array}{ccc}
      -\sin{{\chi}_{G}} & -\cos{{\chi}_{G}} & 0 \\
      \cos{{\chi}_{G}} & -\sin{{\chi}_{G}} & 0\\
      0 & 0 & 0\\
    \end{array}
  \right]
  \dot{\chi}_{G} e
  + R_z({\chi}_{G})
\dot{e}\,.
\end{equation}
We substitute the kinematics of both the vehicle models, equations~\eqref{sys:UAVdyn_Wind} and~\eqref{sys:UGVdyn} in equation~\eqref{eq:noname}, that results in 
\[
\begin{split}
R_z(\chi_A) R_y(\gamma_A) 
\begin{bmatrix}
v_{A}\\ 0\\0
\end{bmatrix}
=&  \,
R_z(\chi_G)
\begin{bmatrix}
v_G\\ 0\\0
\end{bmatrix}
-R_z(\chi_G)
\begin{bmatrix}
e_y\\ -e_x\\ 0
\end{bmatrix} 
v_G \sigma_G + R_{z}(\chi_G) \dot{e} \, ,
\end{split}
\]
that is
\[
R_z(\chi_G)^T R_z(\chi_A) R_y(\gamma_A) 
\begin{bmatrix}
v_{A}\\ 0\\0
\end{bmatrix}
= 
\begin{bmatrix}
(1 - e_y \sigma_G) v_G\\ 
e_x \sigma_G v_G\\
0
\end{bmatrix}
+\dot{e} \, .
\]
Now it is straightforward to compute the expression of $\dot{e}$ as 
\begin{equation} \label{eq:kinem_error}
\dot{e}
  =
 {R}_{z}(\chi_G)^T R_z(\chi_A)R_y(\gamma_A) 
\begin{bmatrix}
v_{A}\\ 0\\0
\end{bmatrix}
-
\begin{bmatrix}
(1 - e_y \sigma_G) v_G\\ 
e_x \sigma_G v_G\\
0
\end{bmatrix}.
\end{equation}
Equation \eqref{eq:kinem_error} describes the kinematic position error of the UAV with respect to the UGV. 
Defining the course angle error as $e_{\chi} = \chi_A - \chi_G$, the speed error as $e_v = v_{A} - v_G$, the flight path error as $e_\gamma = \gamma_A$, and the roll error as $e_\phi = \phi_A$, the coupled nonlinear system \eqref{sys:UAVdyn_Wind},~\eqref{sys:UGVdyn}, can be written with respect to the new set of coordinates 
$(\xb, \ub) = \left(e_x, e_y, e_z, e_v, e_\gamma, e_\chi, e_\phi,  v_G, s_G, u_1, u_2, u_3, u_4 \right)$ as
\begin{equation} \label{sys:CoupledDyn}
    \begin{split}
	\dot{e}_x &= (e_v + v_G) \cos{e_\chi} \cos{e_\gamma} -(1-\sigma_G e_y) v_G ,\\    
	\dot{e}_y &= (e_v + v_G) \sin{e_\chi} \cos{e_\gamma}  -e_x\sigma_G v_G ,\\    
	\dot{e}_z &= -(e_v + v_G) \sin{e_\gamma},\\    
	\dot{e}_v &= \frac{u_1-D}{m}-g\sin e_\gamma - u_4,\\
	\dot{e}_\gamma &= \frac{1}{(e_v + v_G)}\left(\frac{L\cos \phi_A}{m}-g\cos e_\gamma\right),\\
	\dot{e}_\chi &= \frac{1}{(e_v + v_G)\cos e_\gamma}\left(\frac{L\sin\phi_A \cos \chi_c}{m}\right) - \sigma_G v_G, \\
	\dot{e}_\phi &= u_2 , \\
	\dot{v}_G &= u_4 , \\
	\dot{s}_G &= v_G . \\
    \end{split}
\end{equation}

Given the coupled UAV-UGV dynamics~\eqref{sys:CoupledDyn} and the constraints~\eqref{sys:UAVConstraints} and~\eqref{sys:UGVConstraints}, 
we introduce two additional constraints. 
First, the vertical error coordinate, $e_z$, must be semi-negative. 
In other words, we are avoiding UAV collision with the ground (since the UGV altitude is zero). 
Second, for the physical docking at the rendezvous point, we define the rendezvous constraint based on the course angle error. At the rendezvous point (i.e., when the kinematic error components, $e_x$, $e_y$, $e_z$, are zero) the course angle error has to be less than a given tolerance, $\bar{e}_{\chi}$. 
Specifically, the following two constraints are taken into account:
\begin{subequations} \label{sys:DockingConstraints}
\begin{align}
    e_z &\leq 0 \,, \label{sys:DockingConstraints1}\\
    | e_\chi | &\leq \left( \frac{e_x}{\bar{e}_x} \right)^2 + \left( \frac{e_y}{\bar{e}_y} \right)^2 + \left( \frac{e_z}{\bar{e}_z}\right)^2 + \left( \frac{e_\chi}{\bar{e}_{\chi}}\right)^2 \,. \label{sys:DockingConstraints2}
\end{align}
\end{subequations}
Due to the presence of the error coordinates in the right hand side of~\eqref{sys:DockingConstraints2}, if the UAV is far away from the UGV, the course angle error is bounded by a large positive number and the constraint is relaxed. In other words, the rendezvous constraint does not affect the UAV behavior when the UAV and UGV are far away each other. 
This constraint formulation allows us to guide the UAV to the UGV for the successful rendezvous thus avoiding the UAV approaching the UGV in a perpendicular direction.

\subsection{Optimal control problem: a trajectory tracking approach for rendezvous} 
We now formulate the rendezvous problem with respect to the coupled UAV-UGV dynamics~\eqref{sys:CoupledDyn}. 
Motivated by the application scenario depicted in Figure~\ref{fig:scenario2}, we assume that the path of the UGV and the time interval for rendezvous are given. 
Specifically, the UGV can move along a fixed path based on the specific scenario (e.g., a pre-determined area is assigned for docking or landing task). 
The time interval for rendezvous enables the UGV to create a schedule for  service different vehicles operating in the same area. 	For this purpose, the UAV must land onto the UGV between a given time interval $[t_0, T]$.  
Moreover, the UAV is aligned with the UGV at the initial of the rendezvous maneuver (i.e., the longitudinal and lateral position errors are zero at time $t_0$). This initial condition will allow us to predict the time to rendezvous which is an important performance feature of the UAV-UGV trajectory. 
In order to accomplish a successful rendezvous, 
we address the problem of computing rendezvous trajectories by using a nonlinear least squares trajectory optimization technique. 
That is, we consider the following optimal control problem
\begin{equation} \label{pb:OCP_time}
  \begin{split}
    \min_{\xb(\cdot), \ub(\cdot)} & \frac{1}{2} \! \int_{t_0}^{T} \!\!\! \left( \|\xb(\tau) \!- \! \xb^d(\tau)\|_Q^2
    \! + \! \|\ub(\tau)\! - \! \ub^d(\tau)\|_R^2 \right) d\tau  +  \frac{1}{2}\|\xb(T)-\xb^d(T)\|_{P_1}^2\\
    \subj &\; \eqref{sys:CoupledDyn}, \: \texttt{{\small{\emph{dynamics constraints}}}} \\
    &\; \eqref{sys:UAVConstraints}, \eqref{sys:UGVConstraints}, \eqref{sys:DockingConstraints}, \; \texttt{{\small{\emph{state/input constraints}}}}   \\
  \end{split}
\end{equation}
where $(\xb^d(\cdot), \ub^d(\cdot))$ is a desired curve, 
$t_0$ and $T$ are fixed, and 
$Q$, $R$ and $P_1$ are positive definite weighting matrices. %
We address the problem \eqref{pb:OCP_time} numerically by using the projection operator based Newton method for trajectory optimization (\texttt{PRONTO}) with barrier function relaxation, see~\cite{JH:02} and \cite{JH-AS:06} for the details. 
\texttt{PRONTO} is a direct method for solving continuous time optimal control problems. 
It exhibits second order convergence rate to a local minimizer (with a condition on the sufficient closeness of the initial trajectory) satisfying second order sufficient conditions for optimality. 
However, a naive choice of the desired curve and the initial trajectory may lead the algorithm converge to a (local) optimal trajectory that is too far from the desired curve and, therefore, not allow us to perform successful rendezvous between the UAV and the UGV. 
In order to deal with this issue and, at the same time, generate a tunable UGV-UAV trajectory for successful rendezvous, in the next section we design an optimal control based strategy which allows us to effectively solve the optimal control problem \eqref{pb:OCP_time}. 
Notice that a detailed description of \texttt{PRONTO} goes beyond the scope of this paper, while we are interested to show the effectiveness of the rendezvous strategy for the generation of optimal UAV-UGV rendezvous trajectory.

\section{Rendezvous Strategy Based on a Trajectory Optimization Technique}\label{sec:solution}

In this section, we describe the optimal control based strategy for UAV-UGV rendezvous. 
Specifically, we propose a rendezvous strategy based on the following two features: 
i) define a suitable aggressiveness index based on the maximum UAV capability; 
ii) choose a desired state-input curve $(\xb^d, \ub^d)$ based on the decoupled UAV-UGV dynamics. 

\begin{figure}[!ht]
\begin{center}
     \subfloat[$u_{1}(k_{aggr})$.]
{
\begin{tikzpicture}[      
        every node/.style={anchor=south west,inner sep=0pt},
        x=1mm, y=1mm,
      ]   
     \node (fig1) at (0,0)
       {\includegraphics[width=6.5cm]{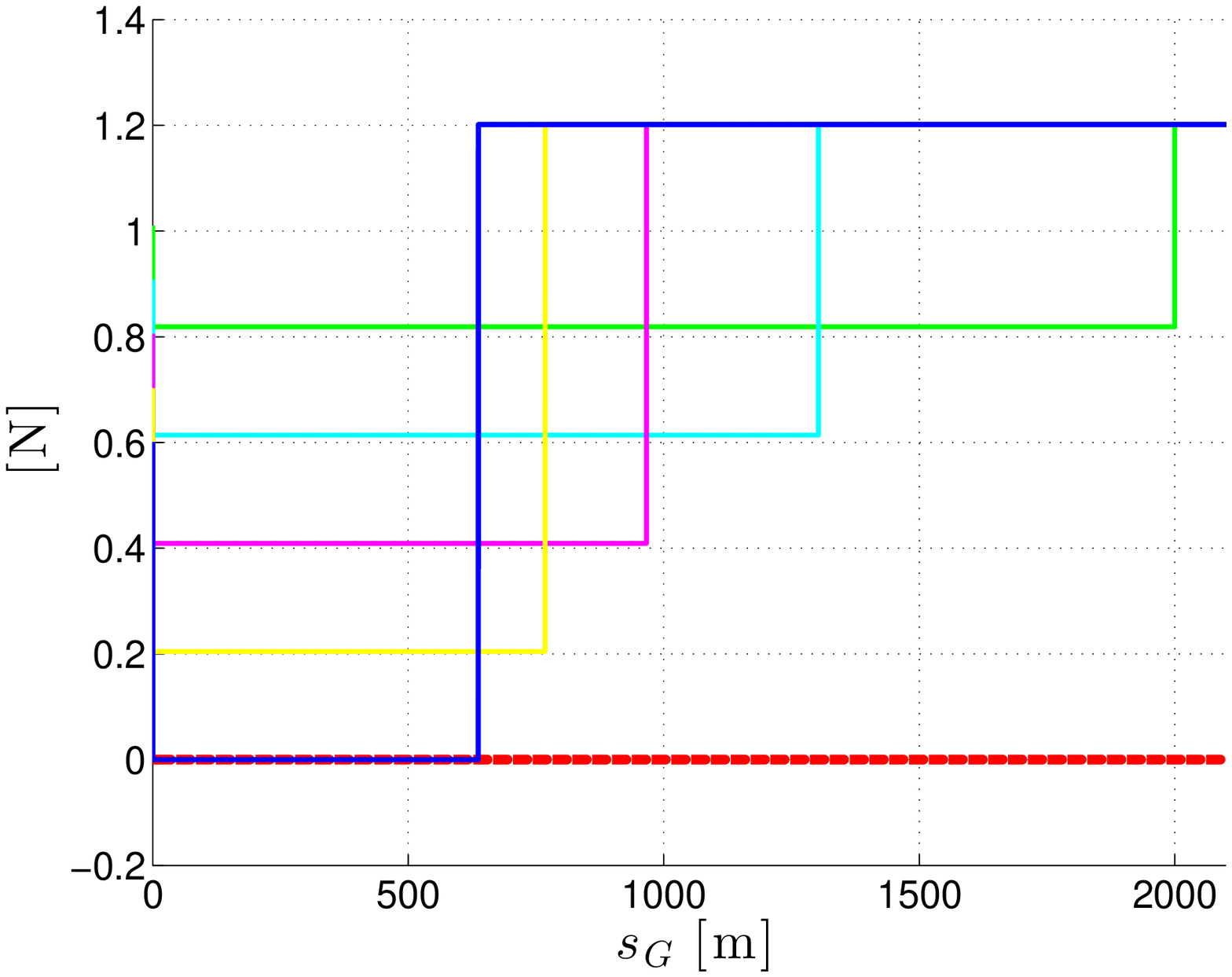} \label{fig:AggrIndex_thrust}};
    \node[draw=none] at (25,12) {\, \small $k_{aggr} = 1$}; 
    \node[draw=none] at (28,17) {\, \small $k_{aggr} = 0.75$}; 
    \node[draw=none] at (34,23) {\, \small $k_{aggr} = 0.5$}; 
    \node[draw=none] at (43,28) {\, \small $k_{aggr} = 0.25$}; 
    \node[draw=none] at (61,34) {\, \small $k_{aggr} = 0$}; 
\end{tikzpicture}  } \,\,\,
     \subfloat[$T^d_{r}(k_{aggr})$.]
     {\includegraphics[width=6.5cm]{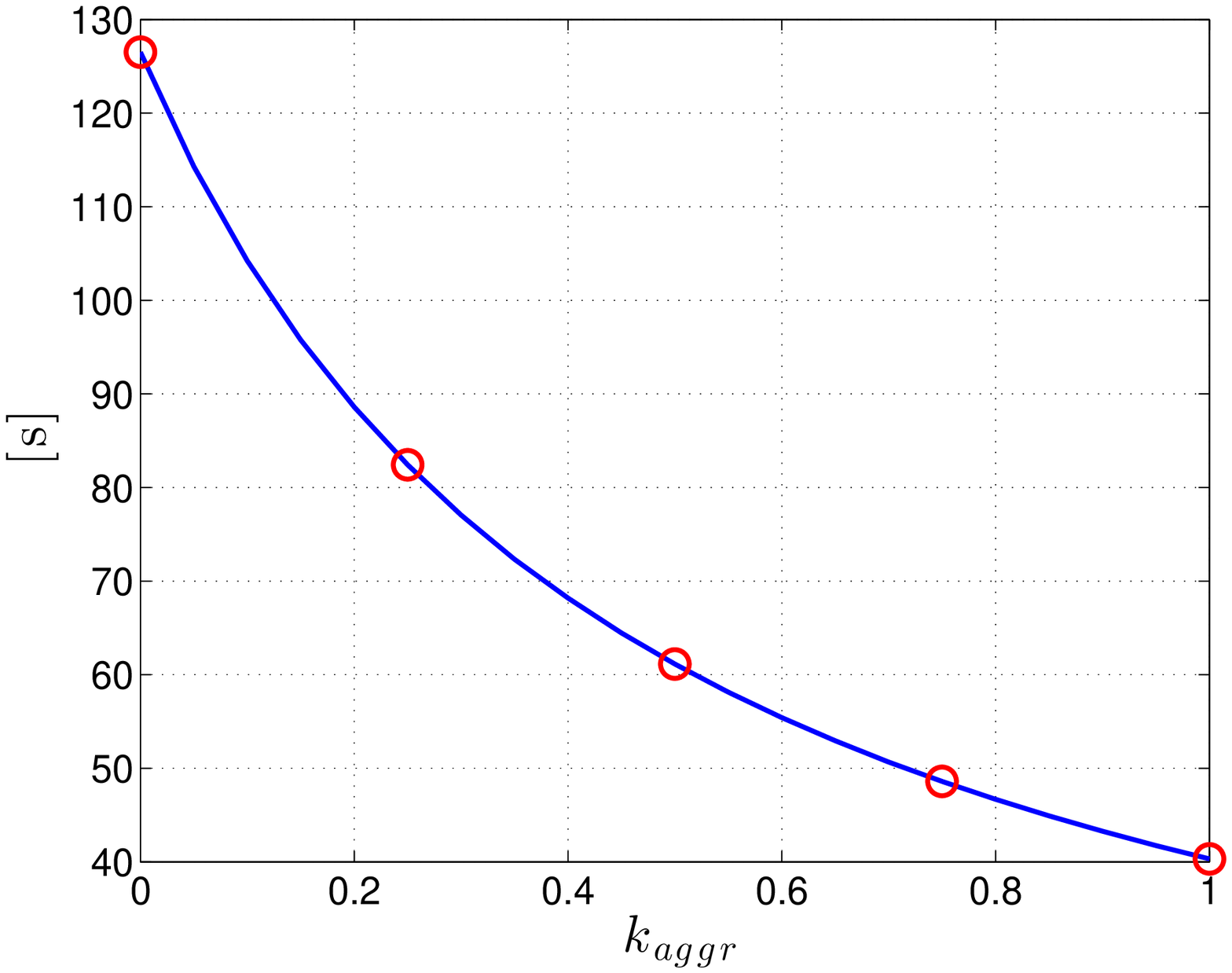} \label{fig:AggrIndex_time}}

    \caption{(a) Thrust and (b) rendezvous time based on the aggressiveness index $k_{aggr} = \{0, 0.25, 0.5, 0.75, 1\}$. 
	}\label{fig:AggrIndex}
\end{center} 
\end{figure}

First, we introduce the aggressiveness index. 
The fixed UGV path is described by the path coordinate $s_G \in [0, s_f]$, where $s_f$ defines the maximum space for the execution of the rendezvous maneuver. 
Let $s_{r}$ be the desired space for the rendezvous maneuver, such that $0 < s_r \leq s_f$. 
Within the rendezvous space interval $s_G \in [0, s_r]$, we set the desired constant flight path angle as 
\begin{equation} \label{eq:gam_des}
{\gamma}_A^d = k_{aggr} \gamma_{1} + (1 - k_{aggr}) \gamma_0 \,, 
\end{equation}
where $k_{aggr} \in [0,1]$ is the aggressiveness index. 
The flight path angle for aggressiveness index equals to zero, i.e., $\gamma_A^d = \gamma_0$, is obtained by imposing the successful execution of the rendezvous maneuver at the maximum space $[0, s_f]$. 
Specifically, the rate of change of  the UAV altitude can be rewritten with respect to $s_G$, i.e., 
${z}'_A = -\sin{{\gamma_A}}$, 
where we use the prime symbol to denote the first derivative of a variable with respect to $s_G$. 
\begin{figure}[htpb]
\begin{center}     
     \subfloat[Lift coefficient vs flight path angle.]
     {\includegraphics[width=6.5cm]{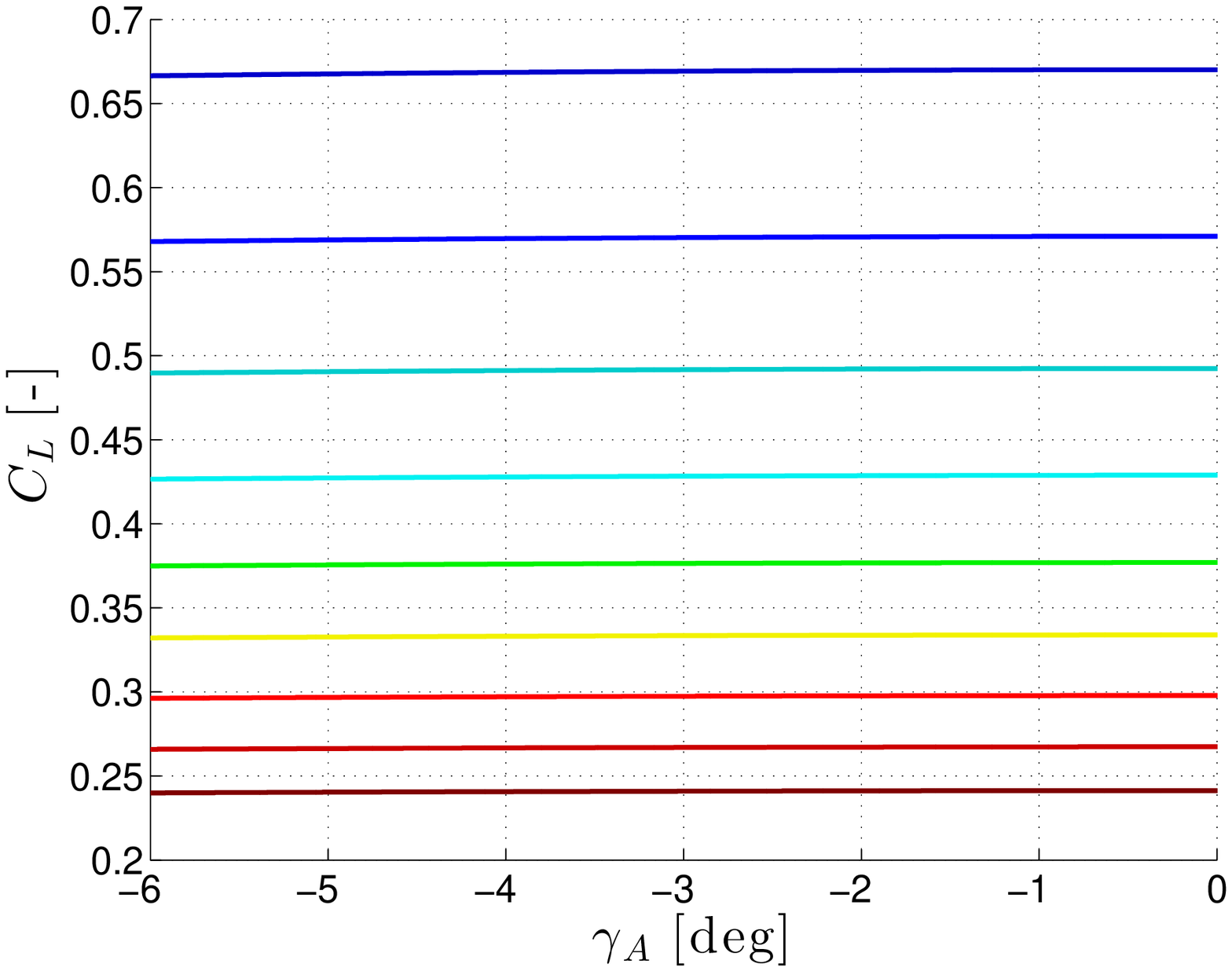} \label{fig:CLtrim}} \,\,\,
     \subfloat[Thrust vs flight path angle.]
{
\begin{tikzpicture}[      
        every node/.style={anchor=south west,inner sep=0pt},
        x=1mm, y=1mm,
      ]   
     \node (fig1) at (0,0)
       {\includegraphics[width=6.5cm]{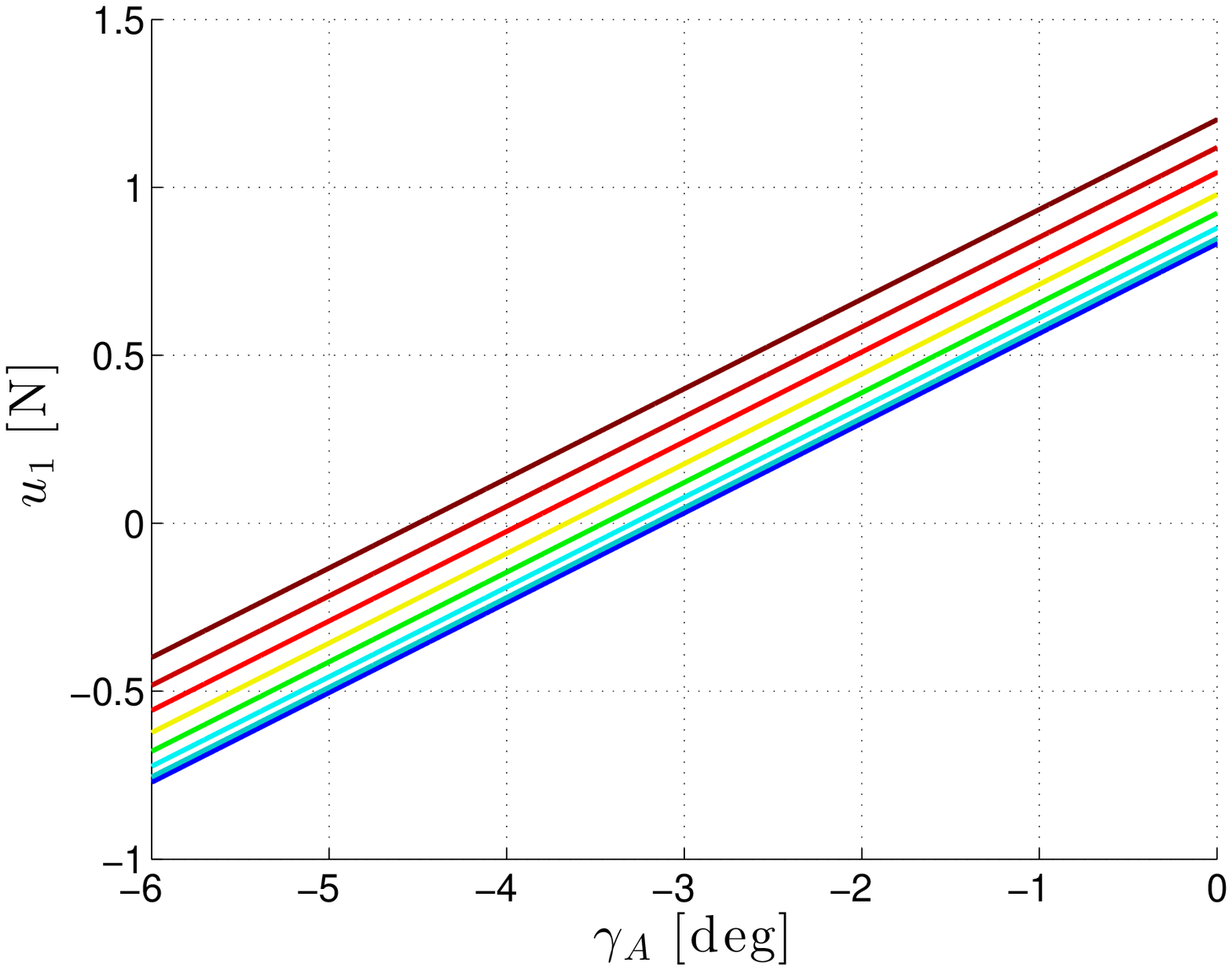} \label{fig:u1trim}};
    \node[draw=none] at (15,42) {\, \small Feasibility region}; 
    \node[draw=none] at (22,38) {\, \small $u_1 \geq 0$}; 
    \draw[red, ultra thick, dashed] (8,23.5) -- (64,23.5);
\end{tikzpicture}  } \,\,\,
     \caption{(a) Trimming trajectories of the fixed-wing UAV with $\phi_A = 0$, $\gamma_A < 0$ and  $v_a = (12, 13, 14, 15, 16, 17, 18, 19, 20)$m/s (blue to red). 
     	(b): negative thrust conditions are avoided. 
     } \label{fig:Trim}
\end{center}
\end{figure} 
By imposing $z_A(s_f) = 0$ and integrating ${z}'_A = -\sin{{\gamma_0}}$, we have 
\[
  \gamma_0 = \arcsin{ \left( \frac{z_0}{s_f} \right)} \,, 
\]
where $z_0$ is the initial UAV altitude at which the rendezvous maneuver begins. 
The flight path angle for aggressiveness index equals to one, i.e., $\gamma_{1}$ in~\eqref{eq:gam_des}, is obtained by analyzing the trimming trajectories of the UAV, i.e., the set of trajectories that can be performed using appropriate constant inputs~\cite{beard2012small}. 
Specifically, we are interested in forward flight with constant descent flight path angle. By setting $\dot{v}_A = \dot{\gamma}_A = 0$, and $\phi_A = 0$ in~\eqref{sys:UAVdyn_Wind}, we have
$C_L  =  \frac{2 m g \cos{\gamma_A}}{\rho S v_a^{2} }$, and $u_1 =  m g \sin{\gamma_A} + \frac{1}{2} \rho S v_a^{2}(C_{D_0} + K_{D/L}C_L^{2})$, see Figure~\ref{fig:Trim}. 
As highlighted in Figure~\ref{fig:u1trim}, the thrust decreases linearly with respect to $\gamma_A$ and becomes negative (and, therefore, unfeasible) for $\gamma_A < -\frac{\rho S v_a^{2}}{2 mg} (C_{D_0} + K_{D/L}C_L^{2})$ (for small values of $\gamma_A$). 
In order to ensure the feasibility of the desired curve, we set 
\[
  \gamma_1 = -\frac{\rho S v_{max}^{2}}{2 mg} \left(C_{D_0} + K_{D/L}\left( \frac{2 m g }{\rho S v_{max}^{2} } \right)^{2}\right) \, .
\] 
In Figure~\ref{fig:AggrIndex_thrust} we show that more close the aggressiveness index is to one, 
more close to the boundary constraint will be the thrust, see Figure~\ref{fig:AggrIndex_thrust}.

Second, we choose a desired state-input curve based on the decoupled UAV-UGV dynamics.  
Exploiting the desired flight path angle based on the aggressiveness index~\eqref{eq:gam_des} and taking into account that the altitude of UGV is constant and equal to zero, the desired vertical error coordinate is given by
\begin{equation} \label{eq:ez_des}
	{e}_z^d(s_G) = z_0 - s_G \sin{(k_{aggr}\gamma_{1} + (1-k_{aggr}) \gamma_0)} , \,\,\, s_G \in [0, s_r] \,.
\end{equation}
For UAV-UGV rendezvous, we have ${e}_z^d(s_r) = 0$ and, therefore, 
\[
s_r = \frac{z_0}{\sin{(k_{aggr}\gamma_{1}+(1- k_{aggr})\gamma_0)}} \,.
\] 
In order to achieve smooth ``docking'', the UAV has to decelerate from the initial speed, $v_0$, to the final speed, $v_f $, with $v_{min}\leq v_f < v_0$. 
To this end, we set the desired speed profile as follows 
\begin{equation} \label{eq:v_des}
{v}^d(s_G) = v_0 + \frac{v_{f} - v_0}{s_r} s_G  \,, s_G\in[0, s_r]. 
\end{equation}
The speed profile $v^d$ is used to time parametrize the path and generate the desired curve for the optimal control problem~\eqref{pb:OCP_time}. 
In particular, given the \emph{space-dependent} desired vertical error~\eqref{eq:ez_des} and speed profile~\eqref{eq:v_des}, the corresponding \emph{time-dependent} desired vertical error and speed profile can be calculated by integrating $dt = ds_G/{v}^d$, i.e.,  
\begin{equation} \label{eq:t_des}
t(s_G) = \int_{0}^{s_G} \frac{ds_G}{{v}^d} \,.
\end{equation}
Now it is straightforward to compute the desired rendezvous time as a function of the aggressiveness index: 
\begin{equation}  \label{eq:Trendezvous}
  \begin{split}
      T_{r}^d 	&= t(s_r) \,,\\
      			&= \frac{s_r}{(v_f - v_0) } \ln{ \frac{v_f}{v_0} } \,. \\
  \end{split} 
\end{equation}
As expected, increasing the aggressiveness index, the desired rendezvous time decreases, see Figure~\ref{fig:AggrIndex_time}. 
It is worth noting that, since the desired speed is strictly greater than zero (note that $v_f > 0$ and $v_0 > 0$), the mapping $s_G \mapsto t(s_G)$ is strictly increasing, so that $t(s_G)$ is well defined. 

Given the desired vertical error, UGV speed profile, and rendezvous time, next we choose the remaining state and input components of the desired curve. 
For successful rendezvous, the desired longitudinal and lateral error coordinates, $(e_x^d, e_y^d)$, course angle error, $e^d_{\chi}$, speed error, $e_{v}^d$, roll error, $e_\phi^d$, and flight path error, $e_\gamma^d$, are set to zero. 
The desired thrust and lift coefficient are chosen by exploiting UAV trim conditions~\cite{beard2012small}. 
In particular, assuming the UAV is in forward flight and constant-altitude flight (i.e., $\gamma_A = 0$, $\phi_A = 0$) with the desired UGV speed profile $v_G^d$, 
and under trim conditions (i.e., $\dot{v}_A = \dot{\gamma}_A = 0$ in~\eqref{sys:UAVdyn_Wind}), we have
\begin{equation} \label{eq:trim}
\begin{split}
u_3^d  &=  \frac{2 m g }{\rho S v_a^{d\,2} }\, ,\\
u_1^d &=  \frac{1}{2} \rho v_a^{d\,2}S(C_{D_0} + K_{D/L}u_3^{d\,2}) \, ,
\end{split}
\end{equation} 
where the desired airspeed, $v_a^d$, is obtained from the desired speed profile and the wind triangle relation~\eqref{eq:vgchigamma_vapsgamma}. 
In order to set the desired thrust and lift coefficient, the wind speed and direction are assumed constant and known to the optimization solver (the wind can be estimated from sensors available in an autopilot module~\cite{langelaan2011wind}).
The desired UAV roll rate, $u_{2}^d$, and the desired UGV acceleration, $u_{4}^d$, are set to zero.  

It is worth noting that, through the definition of the desired vertical error coordinate~\eqref{eq:ez_des}, the rendezvous problem~\eqref{pb:OCP_time} is parametrized by the aggressiveness index. 
The main motivation to use the aggressiveness index is twofold: 
predict the time-to-rendezvous (i.e., equation~\eqref{eq:Trendezvous}) and provide a tool in form of tuning knob (we recall that $k_{aggr} = [0,1]$) which allows one to control the aggressiveness of the UAV trajectory. 

Now, with the desired curve in hand, we design the initial trajectory to initialize \texttt{PRONTO} as follows. 
The UAV is in forward flight, constant-altitude flight equal to $z_0$, and constant speed profile equal to the initial speed $v_0$. 
The UGV is traveling along the pre-determined path with constant speed equal to $v_0$. 
Given the initial UAV and UGV trajectories, 
the initial trajectory for the coupled UAV-UGV dynamics (i.e., it satisfies~\eqref{sys:CoupledDyn}) is obtained by using~\eqref{eq:xtow}. 

We highlight that the desired curve is not a trajectory (it does not satisfied the coupled UAV-UGV dynamics) whereas the initial trajectory is a non-aggressive maneuver, which is easy to compute. This is an important point of the strategy. The desired curve is in fact a guess and we leave \texttt{PRONTO} to take care of the dynamics and state-input constraints and thus compute a trajectory (i.e., satisfying the UAV-UGV dynamics).

\section{Numerical computations} \label{sec:results}

We illustrate the proposed UAV-UGV rendezvous strategy using numerical computations. 
We start with a relatively simple benchmark scenario: the UAV is landing onto the UGV which is moving along a straight line path. 
Then, motivated by the scenario in Figure~\ref{fig:scenario2}, we take into account a $90^\circ$ turn for the UGV path: the strong coupling between longitudinal and lateral dynamics of both UAV and UGV makes the computations particularly challenging and allows us to strengthen the results.
For both scenarios, the rendezvous maneuver starts at $t = 50$sec. The initial ground speed is $18$m/s and the final rendezvous speed is set to $1.15 v_{min}$. We assume planar wind field with wind components $(w_x, w_y, w_z) = (-4.33, 2.5, 0)$. 
It is worth noting that, differently from the approach proposed in~\cite{rucco2016optimal}, we do not tune the $13$ terms in the weighing matrices (they are the same for all the computations). In order to control the aggressivness of the (local) optimal trajectory, we modify only one parameter, i.e., the aggressiveness index $k_{aggr}$.  

\subsection{Rendezvous on a straight line path}
\begin{figure}[htpb]
  \begin{center}
     \includegraphics[width=8cm]{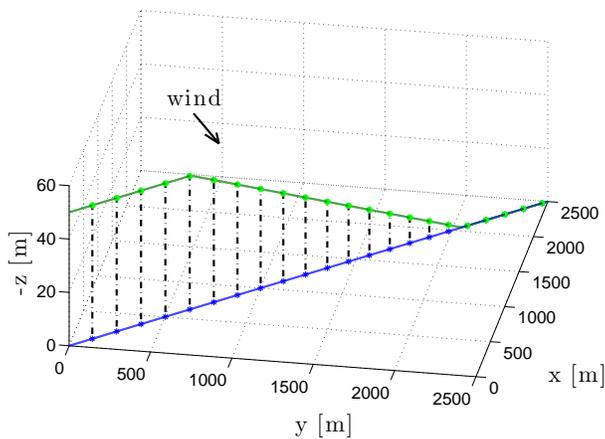}
    \caption{Rendezvous on a straight line path for $k_{aggr} = 0$: (local) optimal path. 
    The blue and green lines represent the UGV and the UAV paths, respectively. 
    } \label{fig:StraightLinePathK0_path}
  \end{center}
\end{figure}

\begin{figure*}[!b]
	\begin{center}
     \subfloat[$-{e}_z$.]
     {\includegraphics[width=5.5cm]{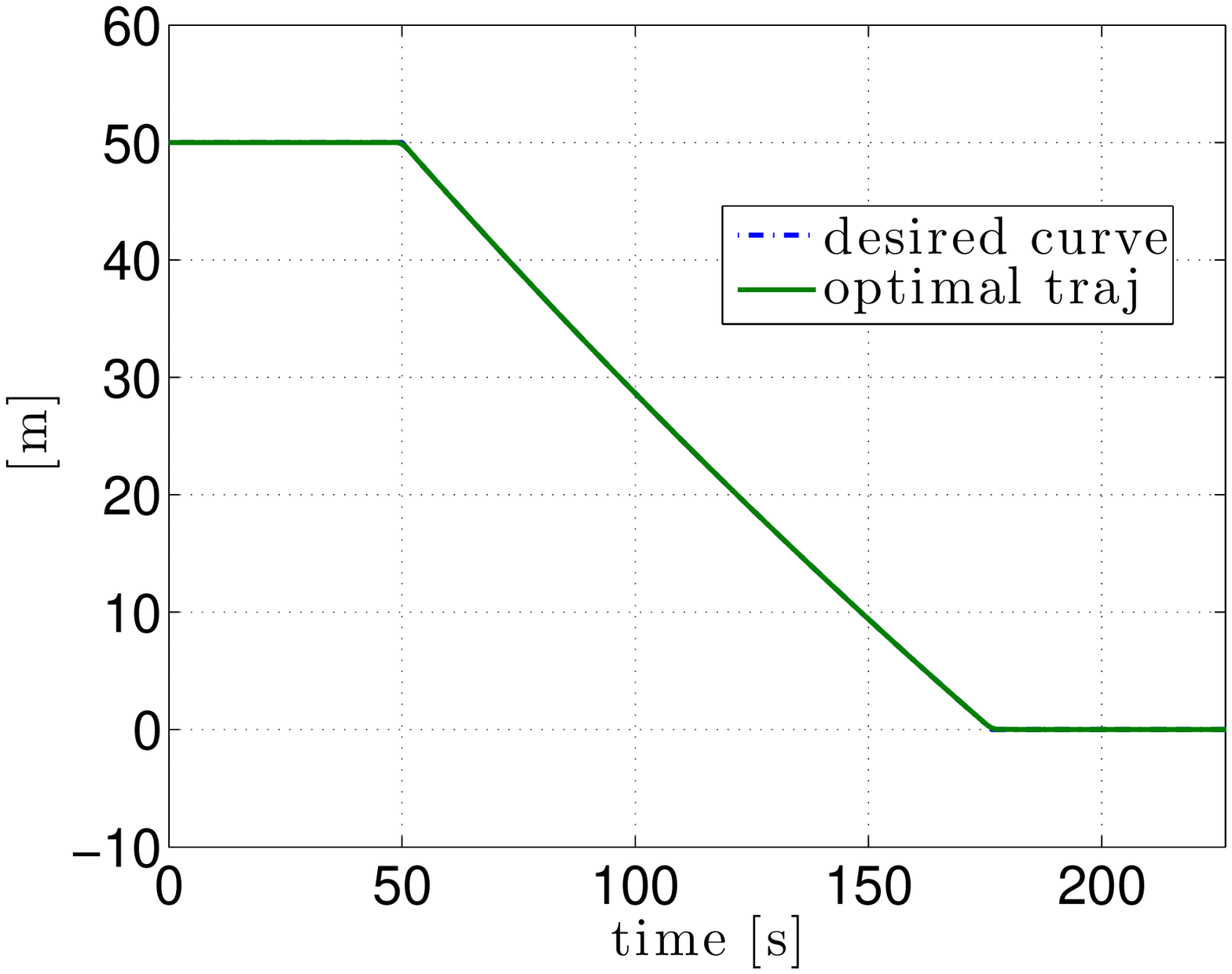} \label{fig:StraightLinePathK0_ez}}
     \subfloat[${e_v}$.]
     {\includegraphics[width=5.5cm]{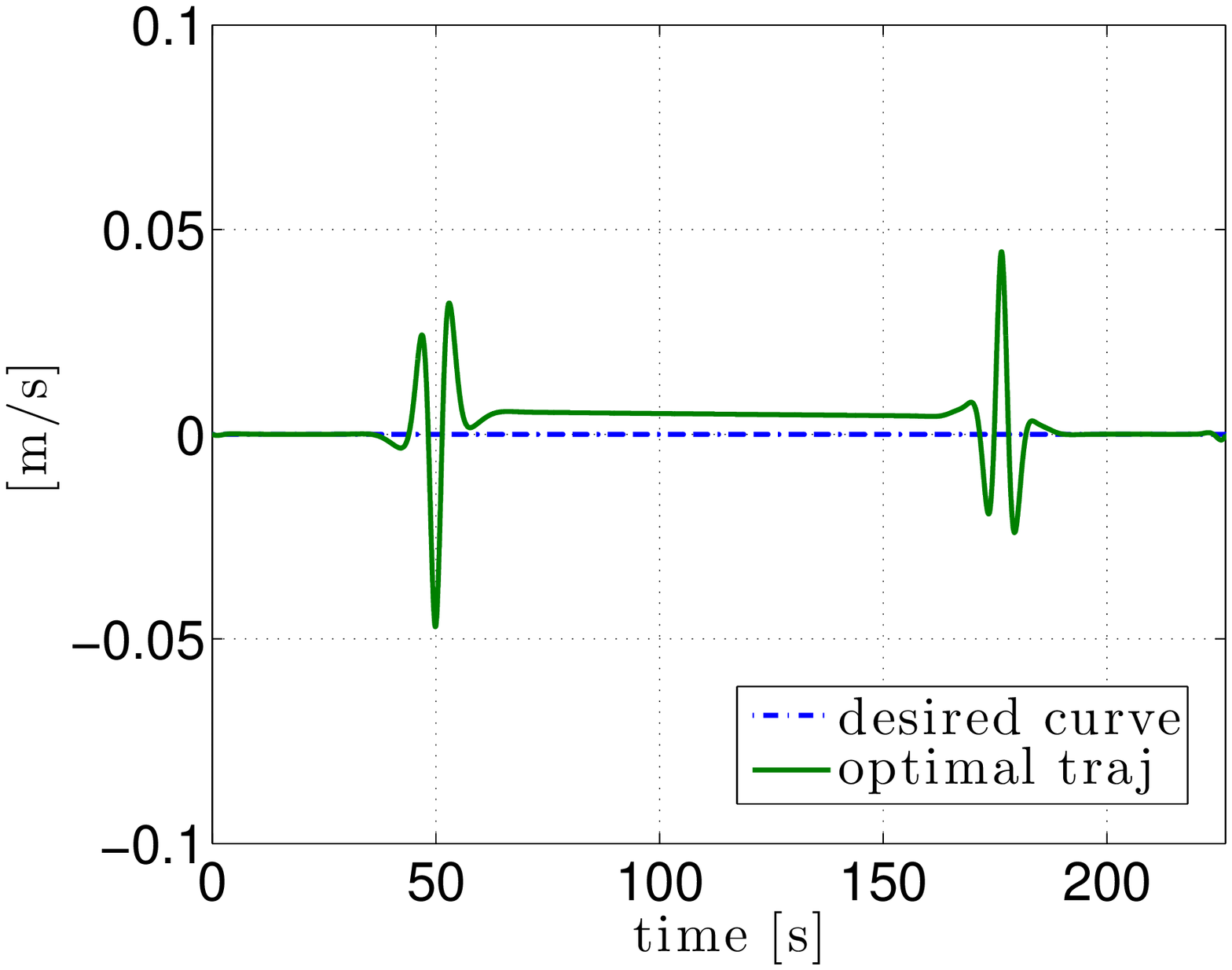}\label{fig:StraightLinePathK0_ev}}
     \subfloat[${e_\gamma}$.]
     {\includegraphics[width=5.5cm]{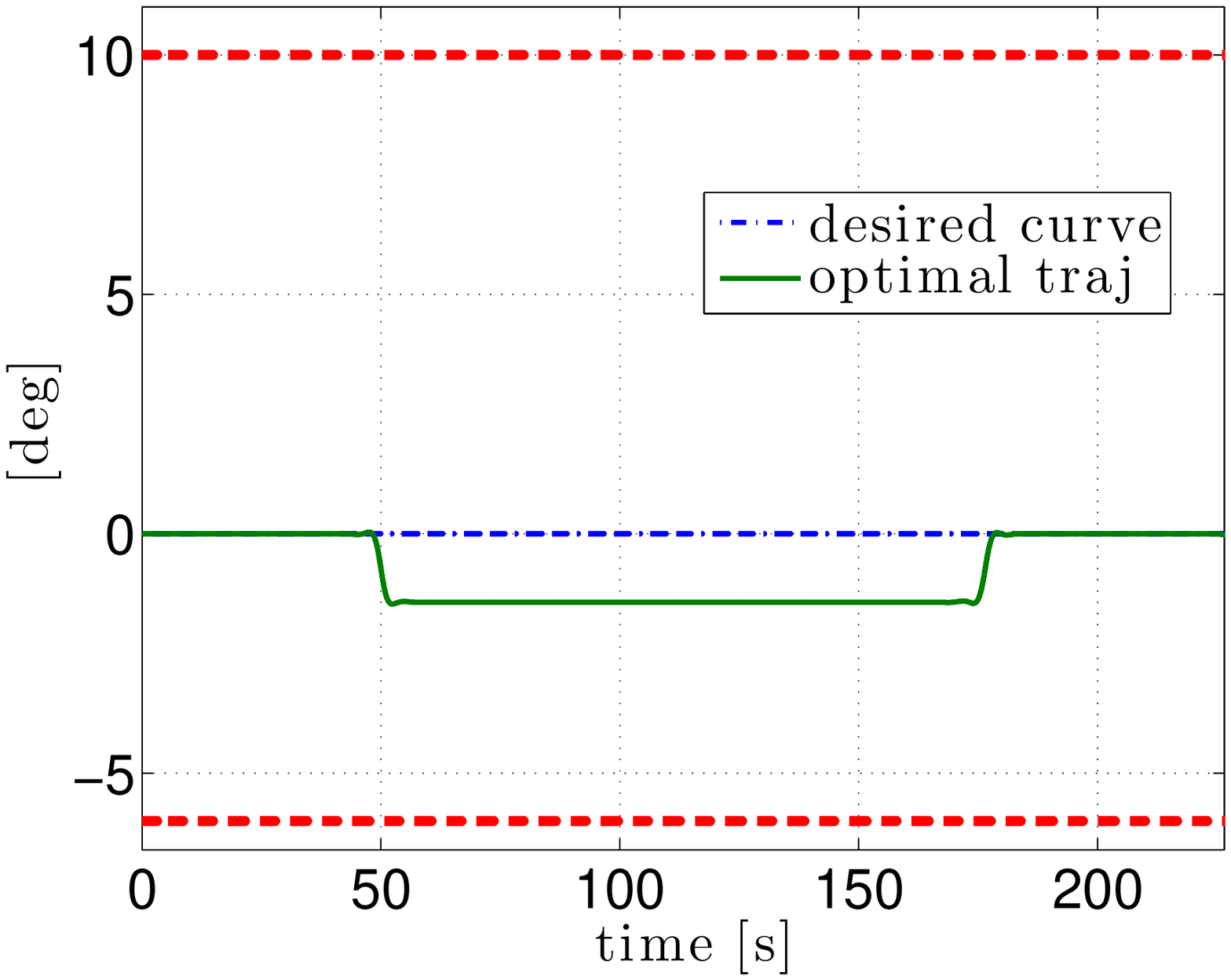} \label{fig:StraightLinePathK0_gam}}

     \subfloat[${u_1}$.]
     {\includegraphics[width=5.5cm]{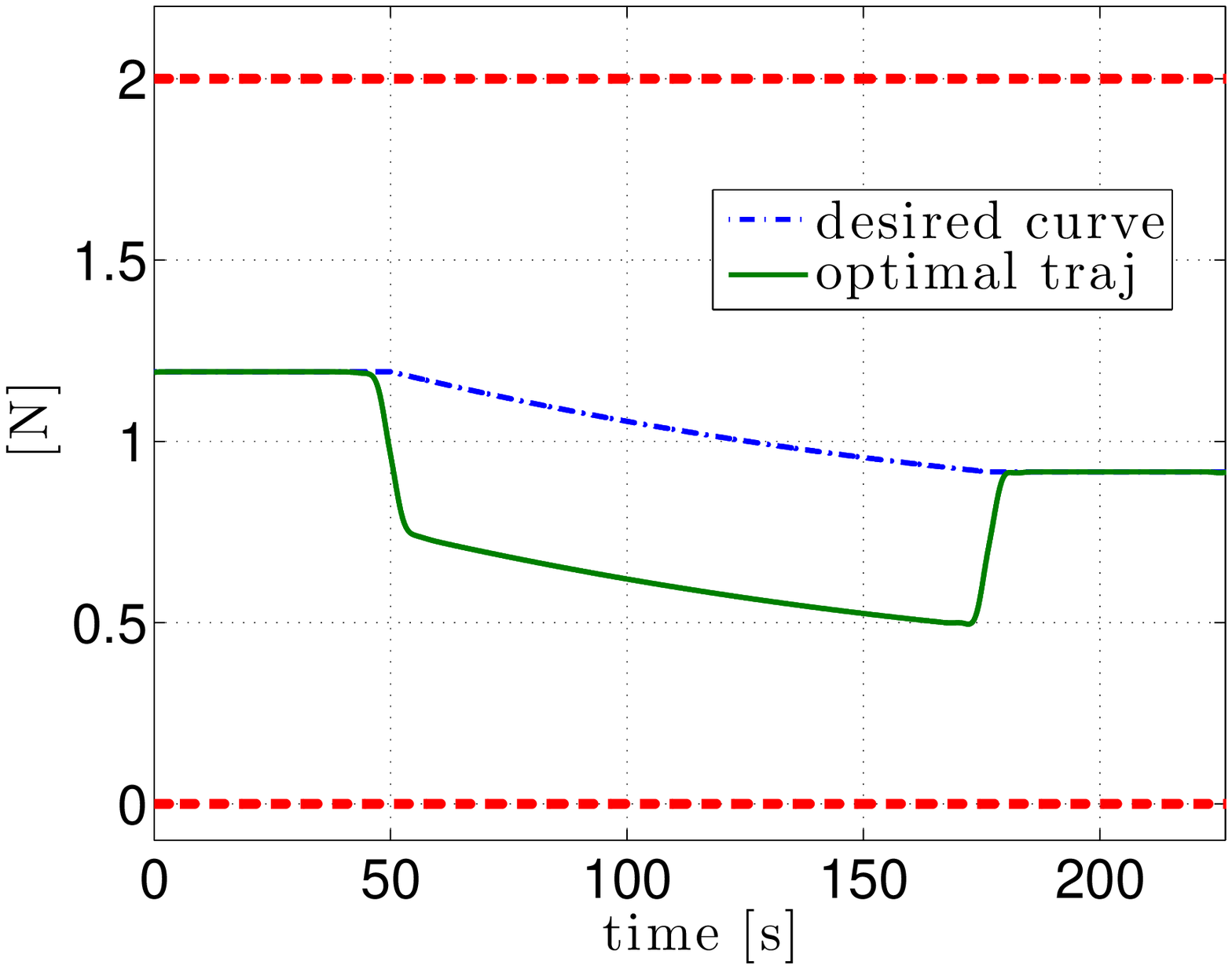} \label{fig:StraightLinePathK0_u1}}
     \subfloat[${u_3}$.]
     {\includegraphics[width=5.5cm]{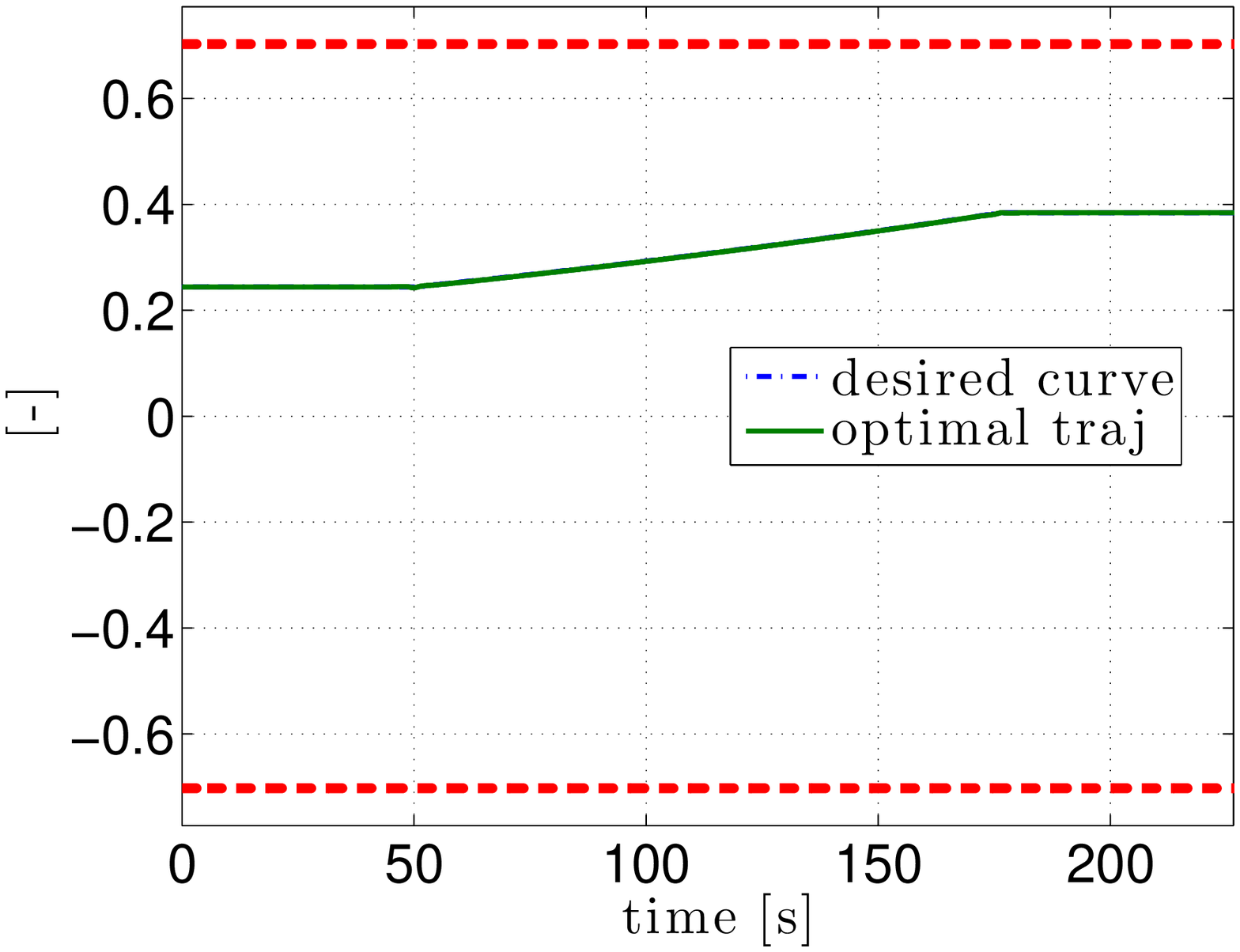} \label{fig:StraightLinePathK0_CL}}
     \subfloat[${n_{lf}}$.]
     {\includegraphics[width=5.5Cm]{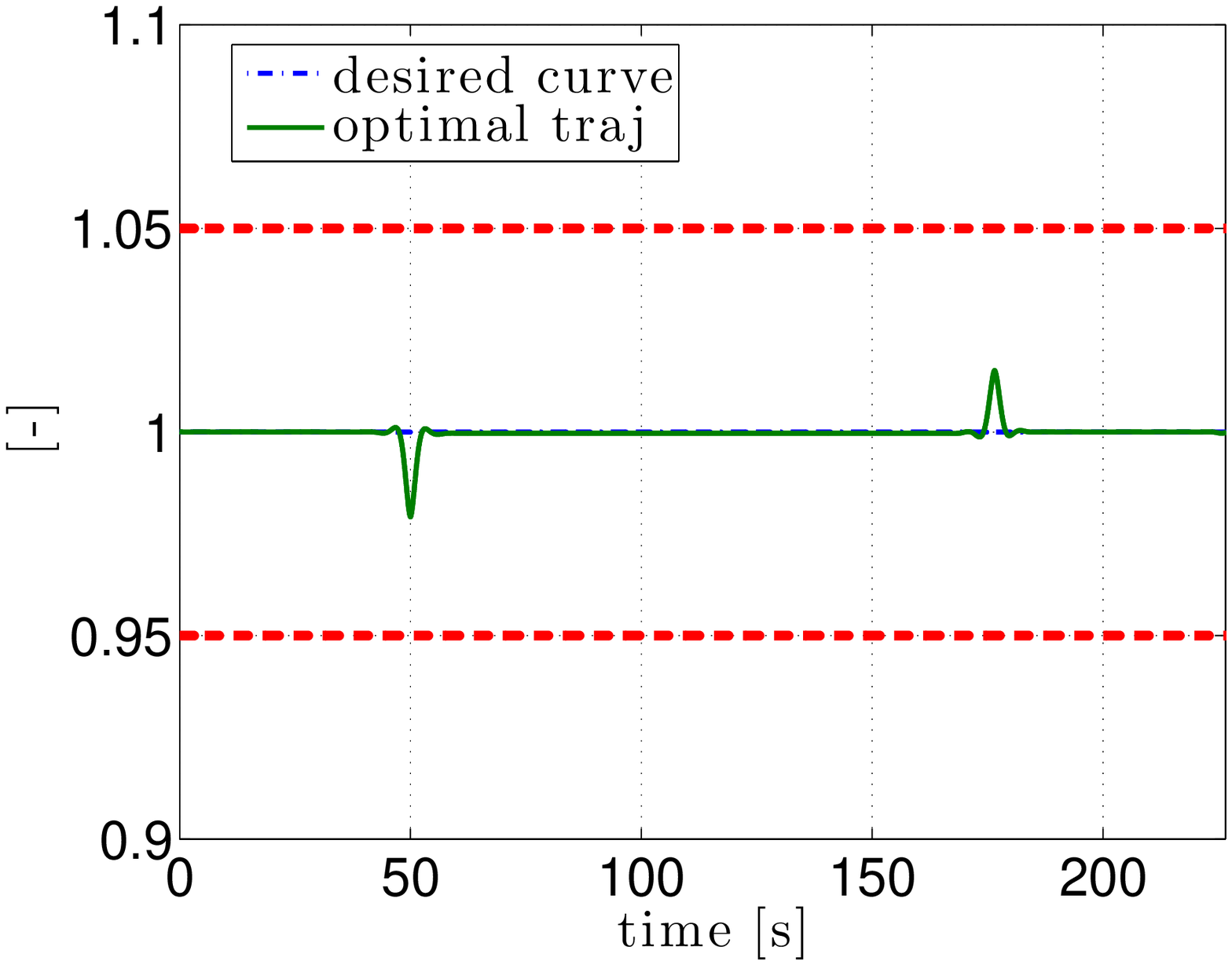} \label{fig:StraightLinePathK0_nlf}}
     \caption{Rendezvous on a straight line path for $k_{aggr} = 0$ (a) vertical error coordinate, (b) error speed, (c) error flight path angle, (d) thrust, (e) coefficient lift, and (f) load factor for $k_{aggr} = 0$. Constraints are in dashed line. } \label{fig:StraightLineK0}
	\end{center}
\end{figure*}

The initial position of the UAV is $(x_A, y_A, z_A) = (0, 0, -50)$, the orientation is $\chi_A = \pi/4$, flight angle and roll angle are $\gamma_A = 0$, $\phi_A = 0$, respectively. 
The initial position and orientation of the UGV are $(x_G, y_G, z_G) = (0, 0, 0)$, and $\chi_G = \pi/4$, respectively. 
The maximum space for the execution of the rendezvous maneuver is $s_f = 2000$m. 
We run \texttt{PRONTO} based on the rendezvous trajectory generation strategy for aggressiveness index equals to $k_{aggr} = 0$. 
We obtain a quadratic convergence rate in the neighbourhood of the solution (we recall that the \texttt{PRONTO} has a structure of a standard Newton method \cite{JH:02}) at each iteration of the algorithm. 
The local optimal trajectory is shown in Figures~\ref{fig:StraightLinePathK0_path} and~\ref{fig:StraightLineK0}.

\begin{center}
\begin{figure*}[!ht]
     \subfloat[$k_{aggr} = 0$.]
     {\includegraphics[width=5.5cm]{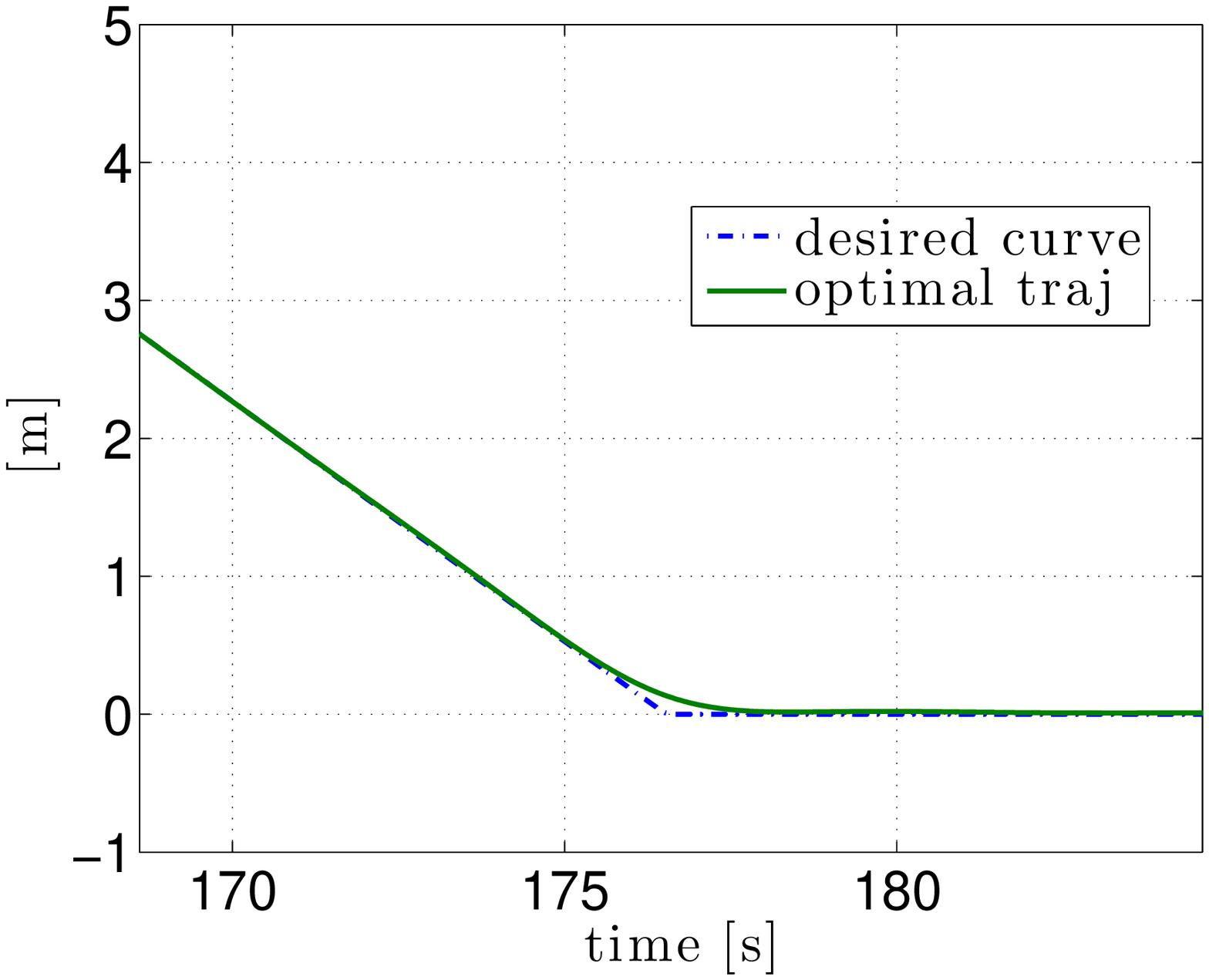} \label{fig:StraightLineZcomparisonK000}}
     \subfloat[$k_{aggr} = 0.25$.]
     {\includegraphics[width=5.5cm]{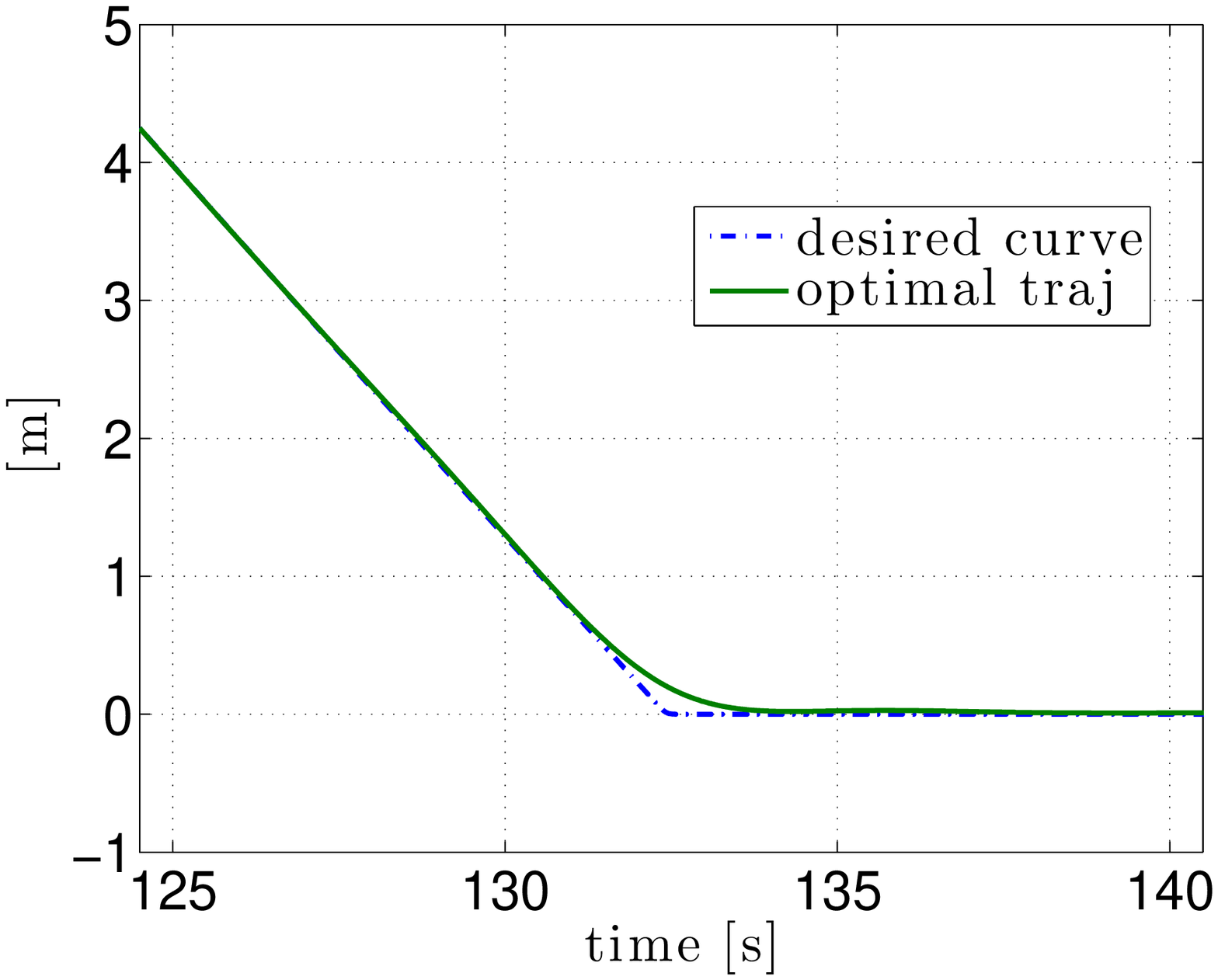}\label{fig:StraightLineZcomparisonK025}}
     \subfloat[$k_{aggr} = 0.5$.]
     {\includegraphics[width=5.5cm]{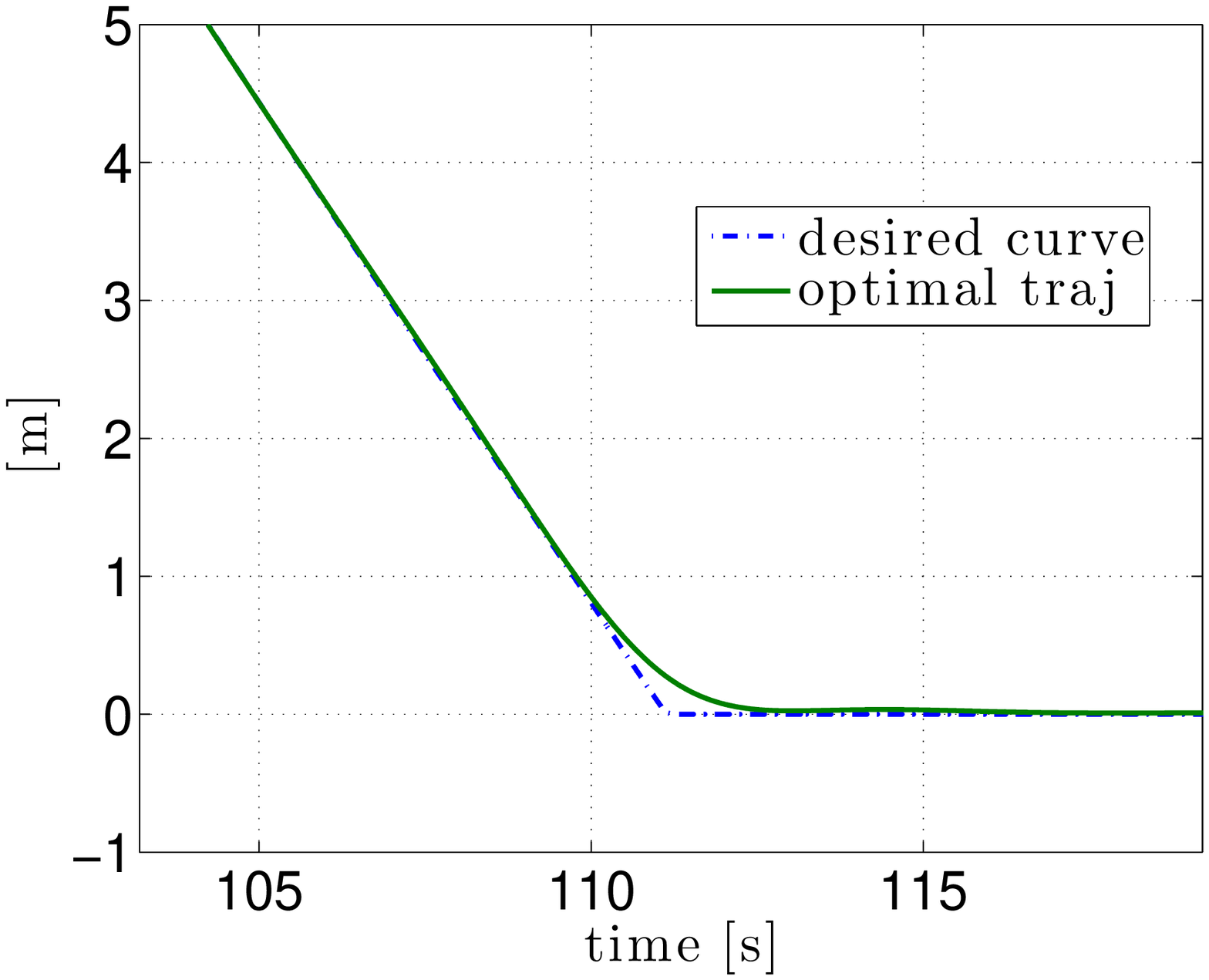} \label{fig:StraightLineZcomparisonK050}}

\begin{center}
     \subfloat[$k_{aggr} = 0.75$.]
     {\includegraphics[width=5.5cm]{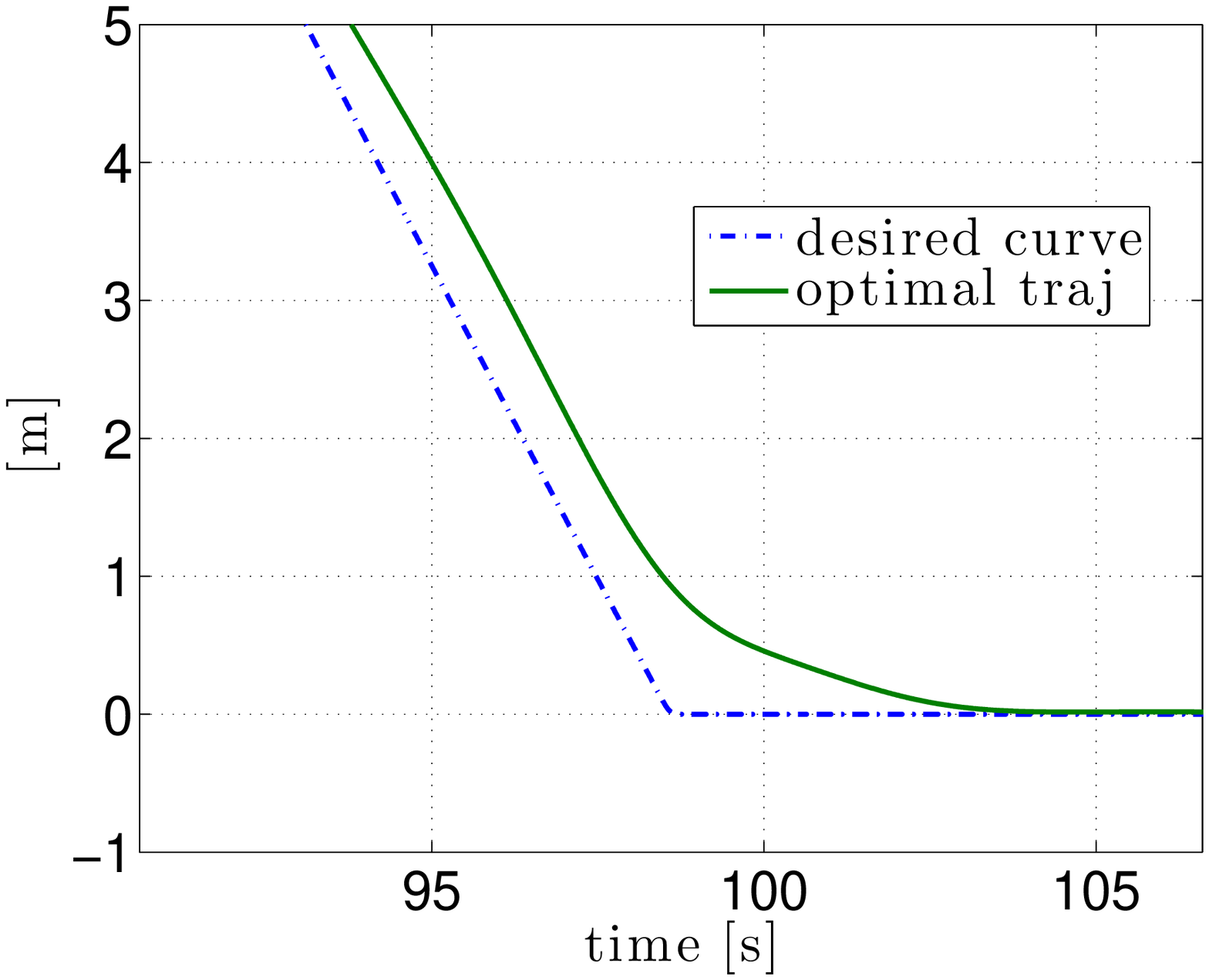} \label{fig:StraightLineZcomparisonK075}}
     \subfloat[$k_{aggr} = 1$.]
     {\includegraphics[width=5.5cm]{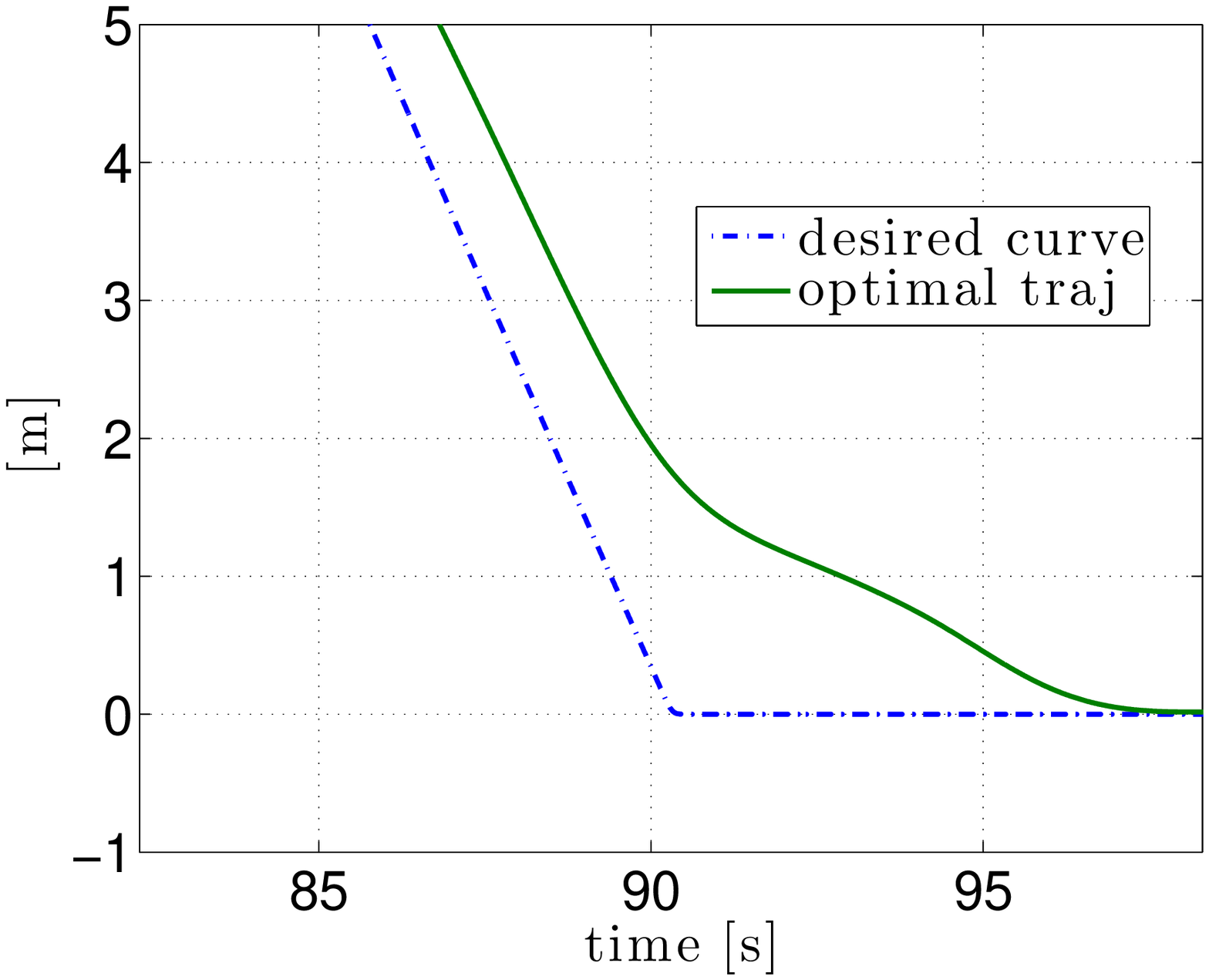} \label{fig:StraightLineZcomparisonK100}}
\end{center}
     \caption{The vertical error coordinate when the UAV is approaching the UGV for different $k_{aggar}$. 
     } \label{fig:StraightLineZcomparisonK}
\end{figure*} 
\end{center}

In Figure~\ref{fig:StraightLinePathK0_path} we show the (local) optimal 3D path traversed by the UAV to rendezvous with the UGV. 
The (local) optimal UAV path is soft: the UAV height is reduced gradually. 
This soft feature is also evident from the trajectory shown in Figure~\ref{fig:StraightLineK0}. 
In fact, the vertical error coordinate and the flight path angle vary smoothly, Figures~\ref{fig:StraightLinePathK0_ez} and~\ref{fig:StraightLinePathK0_gam}, 
and the constraints on thrust, flight path angle, coefficient lift, and normal load are never active, Figures~\ref{fig:StraightLinePathK0_u1},~\ref{fig:StraightLinePathK0_gam},~\ref{fig:StraightLinePathK0_CL},~\ref{fig:StraightLinePathK0_nlf}. 
We observe that the (local) optimal thrust is different from the desired one, see Figure~\ref{fig:StraightLinePathK0_u1}. 
Such a difference is due to the fact that the desired curve is based on trim conditions (i.e., speed transition in forward flight and constant altitude, see~\eqref{eq:trim}) and, thus, does not take into account the change in the flight altitude as well as important dynamic features. 
Finally, we highlight that the rendezvous time is 
$126.7$sec (note that $e_z = -0.1$m for $t = 176.7$sec, see Figure~\ref{fig:StraightLineZcomparisonK000}) 
and the desired rendezvous time is $T_r^d = 126.5$sec, see~\eqref{eq:Trendezvous}.  

\begin{figure*}[!ht]
\begin{center}
     \subfloat[$-{e}_z$.]
     {\includegraphics[width=5.5cm]{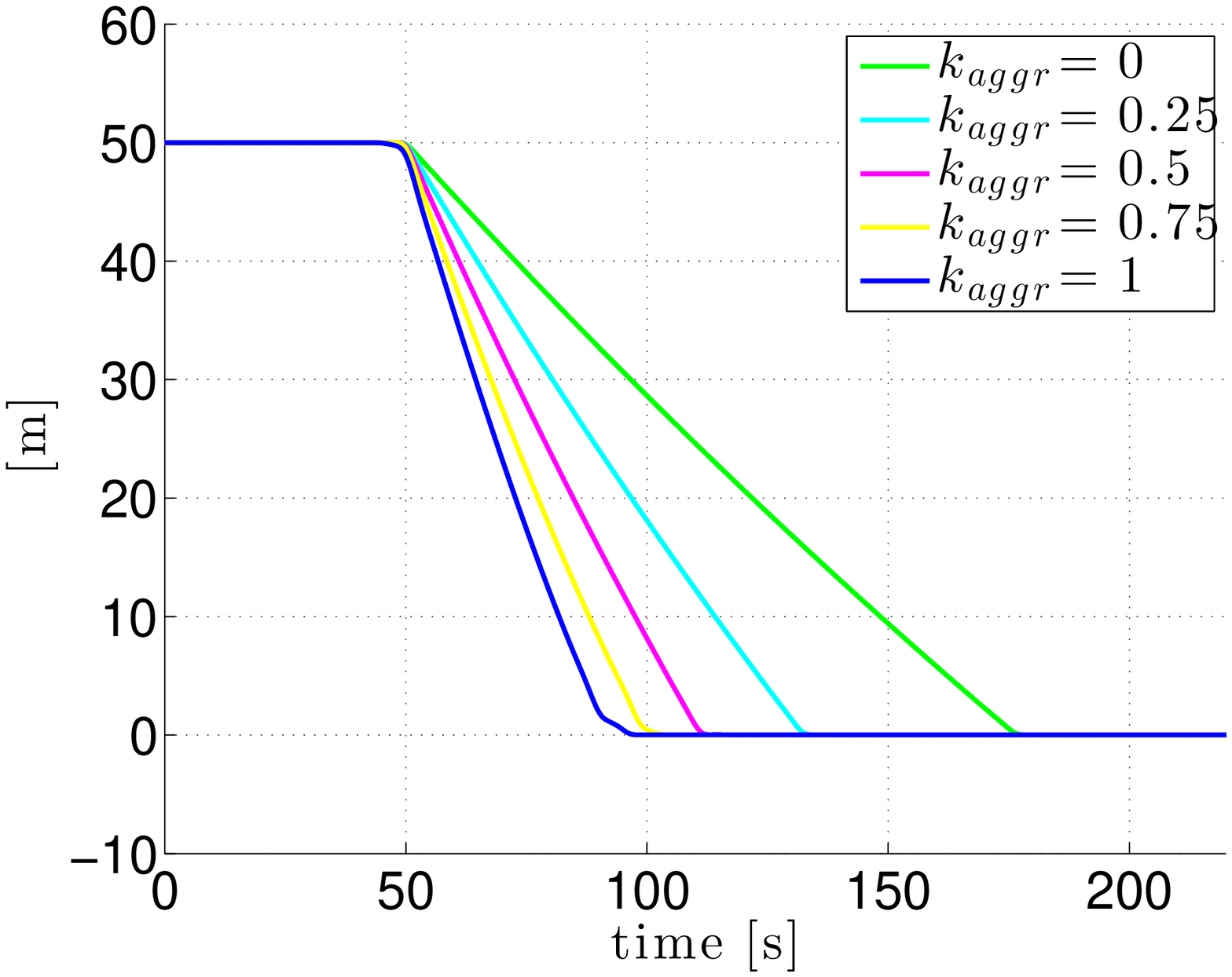} \label{fig:StraightLinePathcomparisonK_ez}}
     \subfloat[${e_v}$.]
     {\includegraphics[width=5.5cm]{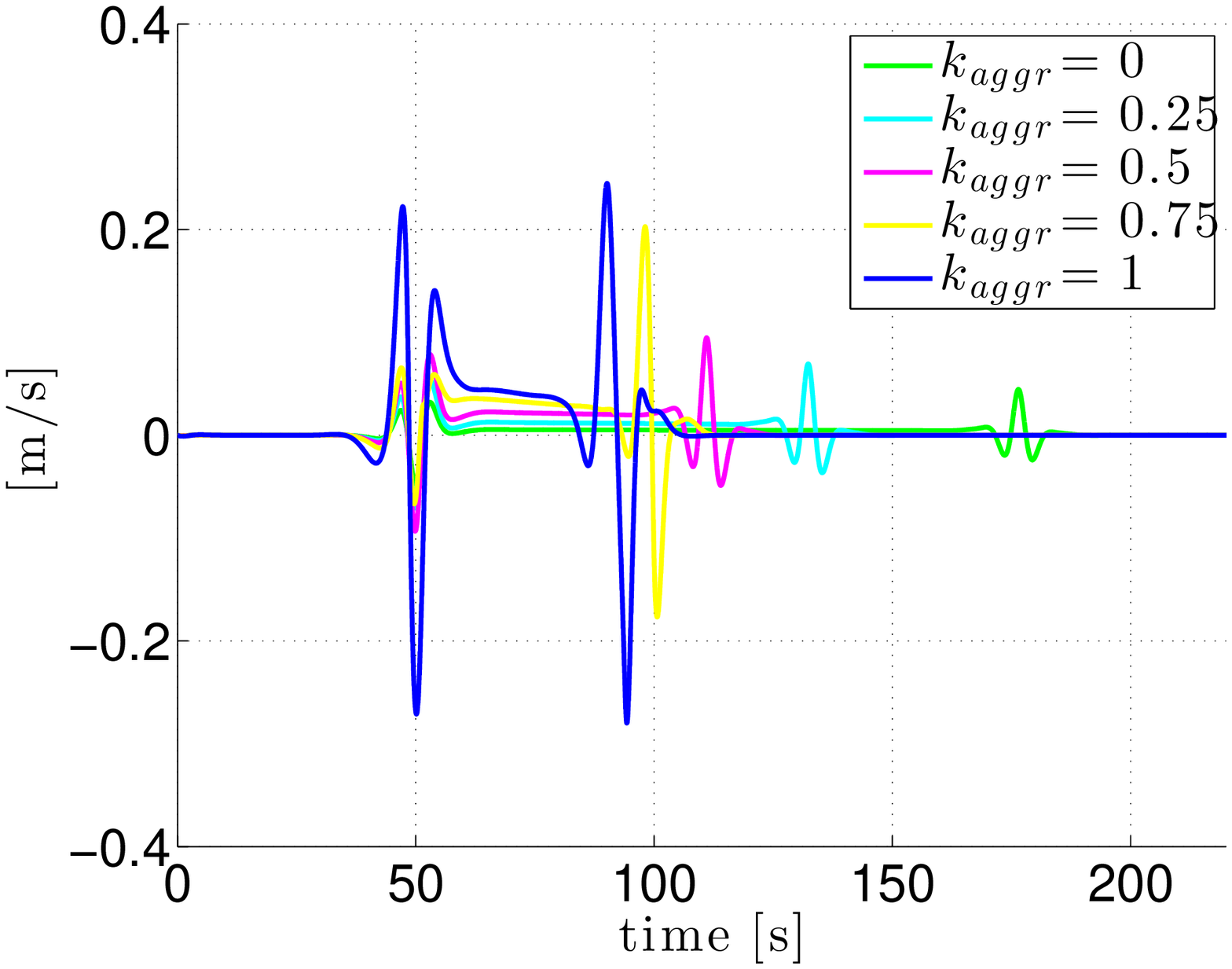}\label{fig:StraightLinePathcomparisonK_ev}}
     \subfloat[${e_\gamma}$.]
     {\includegraphics[width=5.5cm]{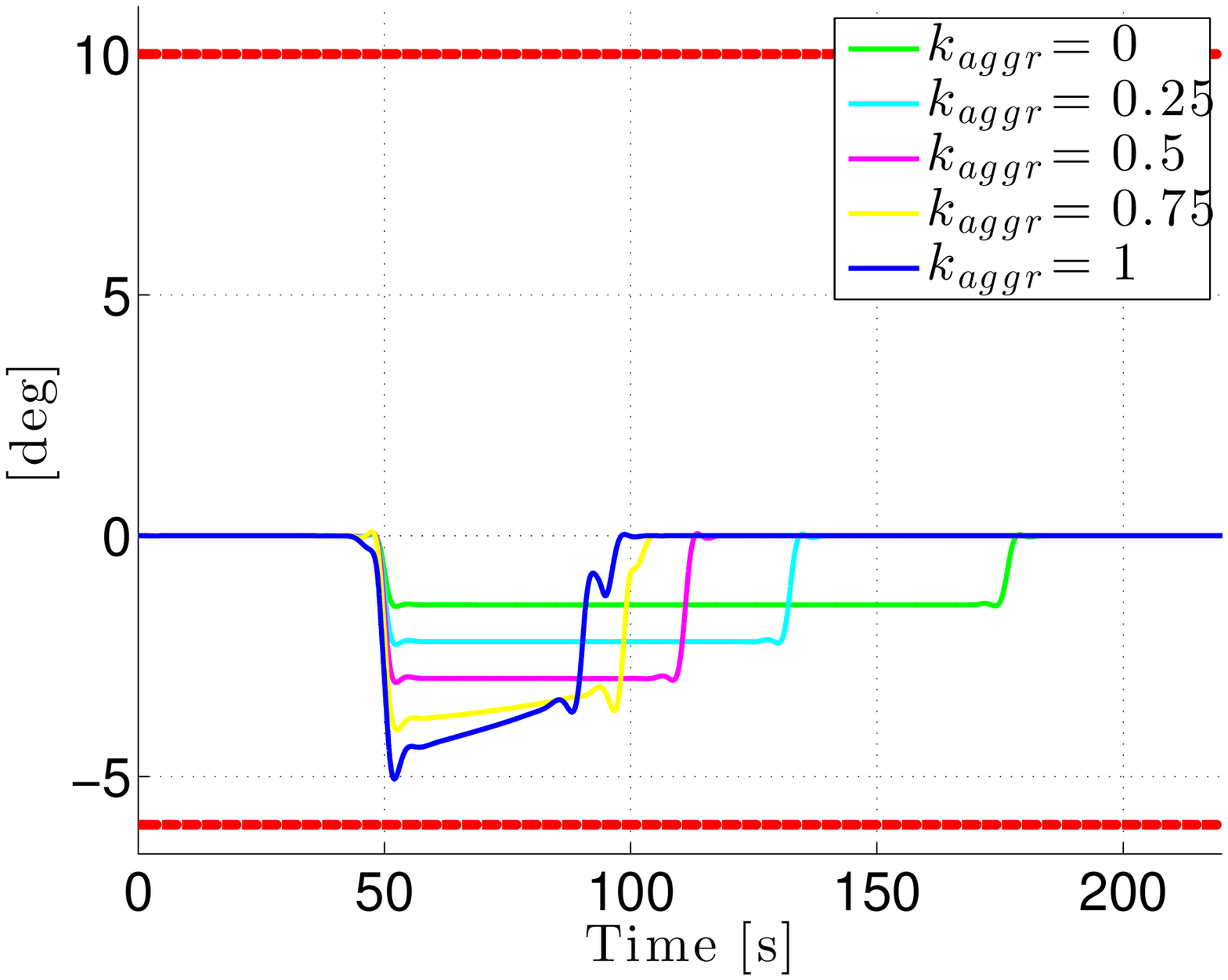} \label{fig:StraightLinePathcomparisonK_gam}}

     \subfloat[${u_1}$.]
     {\includegraphics[width=5.5cm]{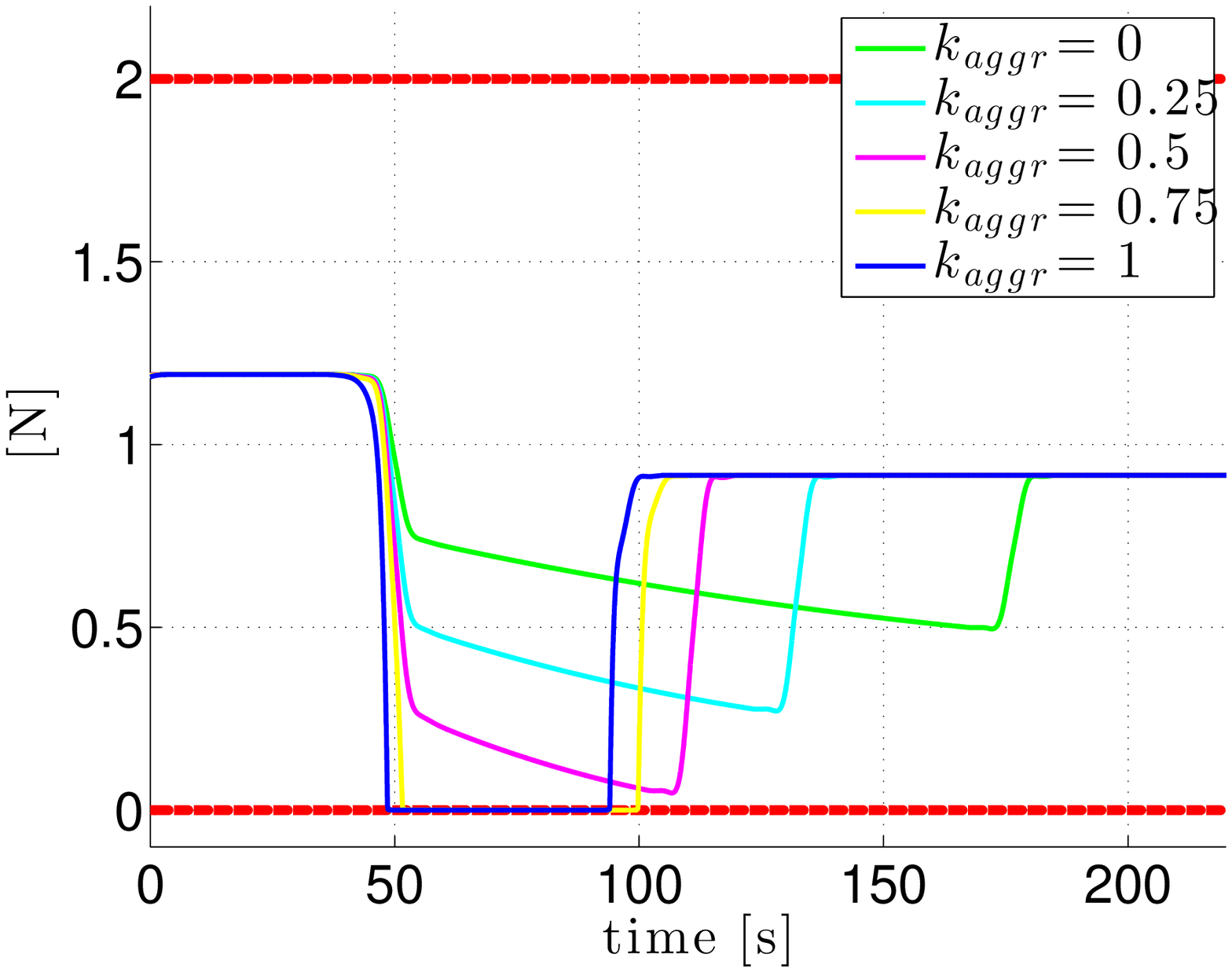} \label{fig:StraightLinePathcomparisonK_u1}}
     \subfloat[${u_3}$.]
     {\includegraphics[width=5.5cm]{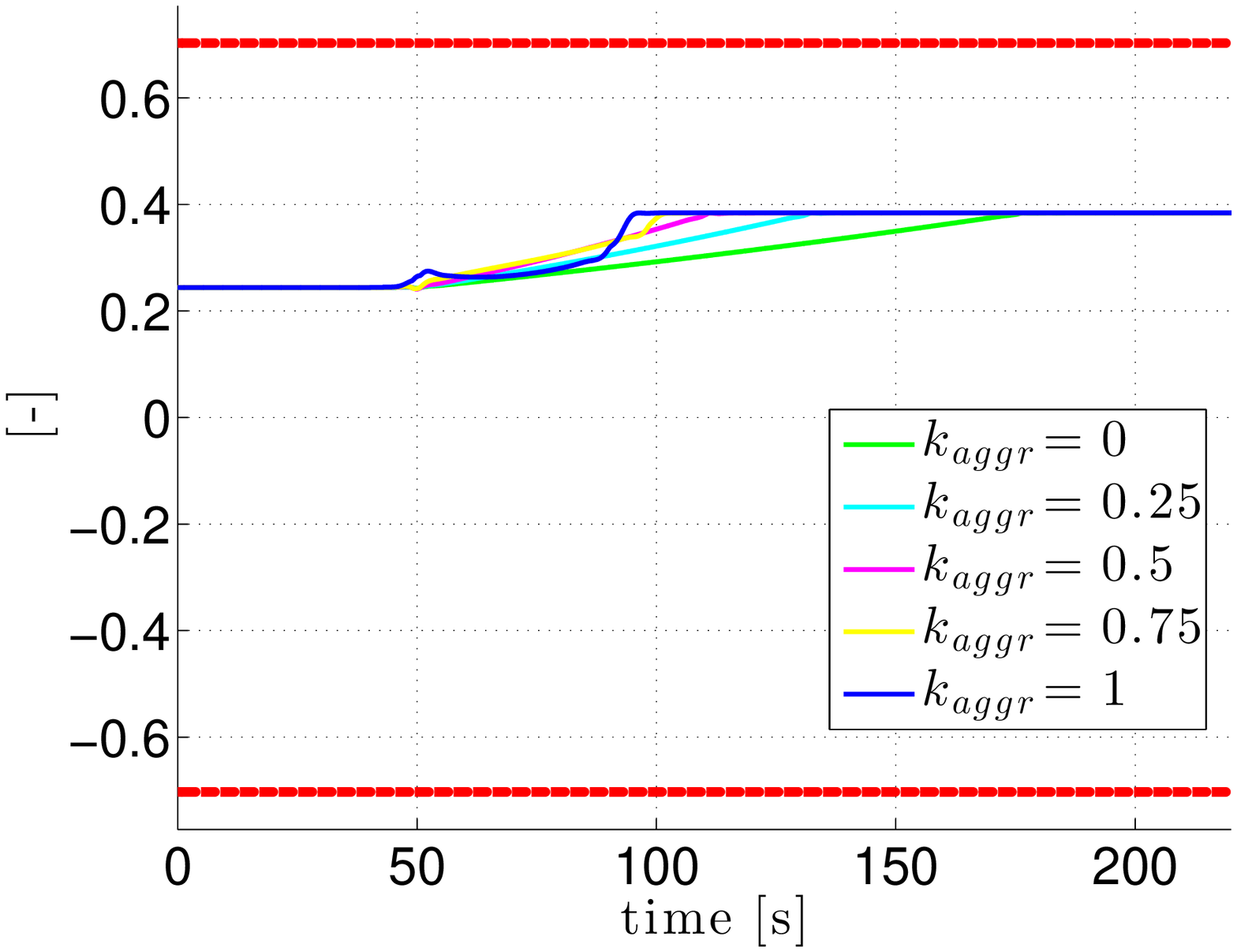} \label{fig:StraightLinePathcomparisonK_CL}}
     \subfloat[$n_{lf}$.]
     {\includegraphics[width=5.5cm]{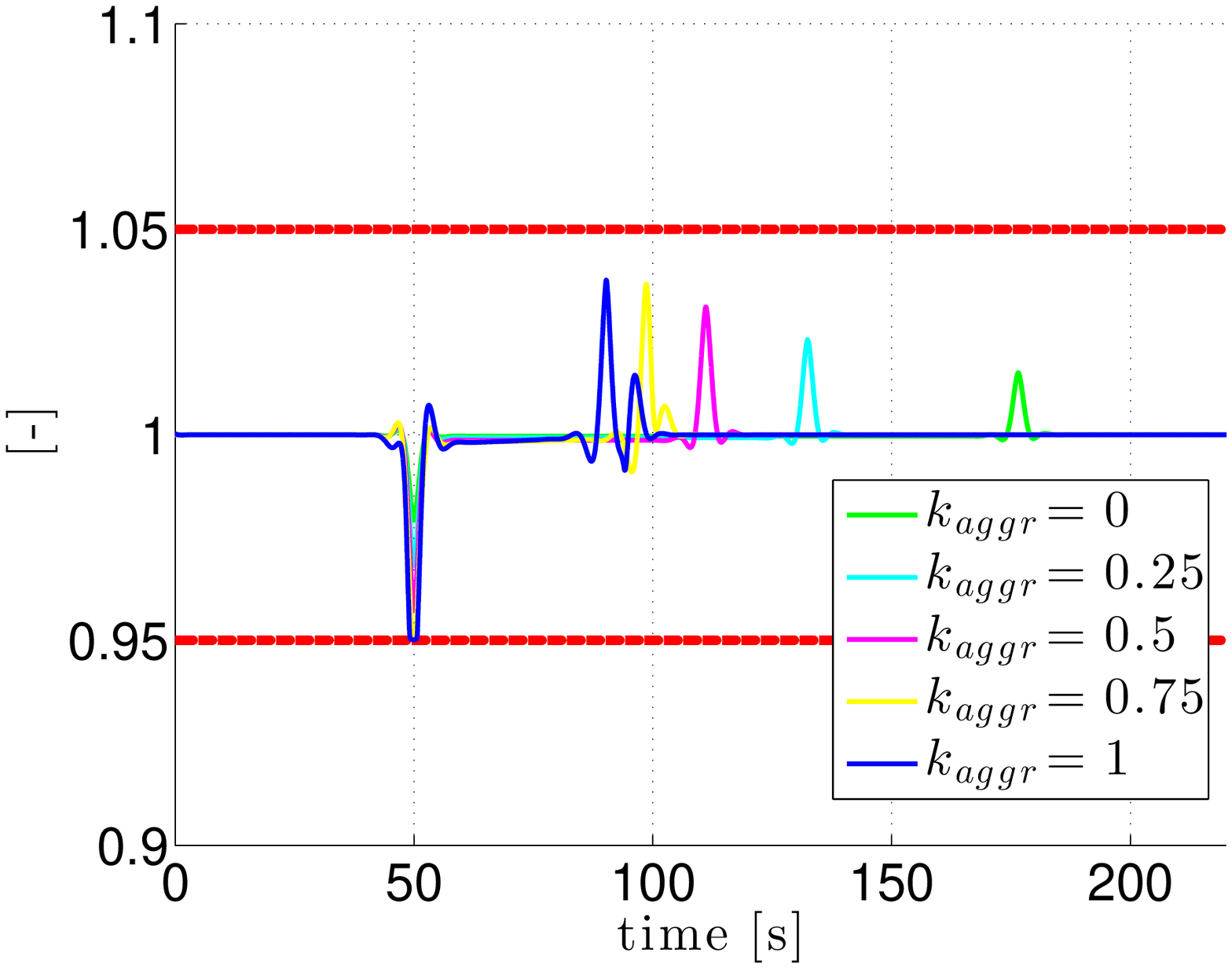} \label{fig:StraightLinePathcomparisonK_nlf}}
     \caption{Rendezvous on a straight line path (a) vertical error coordinate, 
     (b) error speed, 
     (c) error flight path angle, (d) thrust, (e) coefficient lift, and (f) load factor for $k_{aggr} = \{0, 0.25, 0.5, 0.75,1\}$. Constraints are in dashed line. 
     } \label{fig:StraightLinePathcomparisonK}
\end{center}
\end{figure*}

Next, we run the rendezvous trajectory generation strategy for aggressiveness index equals to $k_{aggr} = \{0.25, 0.5, 0.75,1\}$ and we compare the (local) optimal trajectories in Figure~\ref{fig:StraightLinePathcomparisonK} (for the sake of completeness, we include the rendezvous trajectory obtained with $k_{aggr} = 0$). 
For $k_{aggr} = 1$ the (local) optimal UAV path (blue line in Figure~\ref{fig:StraightLinePathcomparisonK_ez}) is aggressive. 
By aggressive, we mean that the several constraints are active during the rendezvous maneuver. 
Indeed, the thrust is zero for almost all the rendezvous maneuver, Figure~\ref{fig:StraightLinePathcomparisonK_u1},
and the constraint on the normal load is active at the beginning of the maneuver, Figure~\ref{fig:StraightLinePathcomparisonK_nlf}. 
Moreover, we observe a (relative) high speed error during the maneuver, Figure~\ref{fig:StraightLinePathcomparisonK_ev}. Such a difference in the ground speed between the UAV and the UGV is due to the wind (which affects only the UAV) and the fact that the airspeed is constrained. 
Next, we highlight two interesting features of the (local) optimal trajectory for $k_{aggr} = 1$. 
First, at the beginning of the rendezvous maneuver, the UAV decreases the thrust and, at the same time, increases the lift coefficient (see the kink at about $t = 50$ sec in Figure~\ref{fig:StraightLinePathcomparisonK_CL}, blue line). 
In this way, the airspeed decreases at about $t = 50$sec and, immediately after, increases thus reaching its maximum value. Such aggressive maneuver allows the UAV to take a steep dive towards the UGV as shown in Figure~\ref{fig:StraightLinePathcomparisonK_ez}. 
Second, once the UAV ground-speed reaches the desired value of $13.8~m/s$ (i.e., $1.15 v_{min}$), the UAV needs to maintain this speed and hence it requires thrust which is increased from zero to the desired value as shown in the Figures~\ref{fig:StraightLinePathcomparisonK_u1}. Similar behavior is observed in the lift coefficient, Figure~\ref{fig:StraightLinePathcomparisonK_CL}. 
Such a sharp variation in the thrust and lift coefficient highlights an important transient behavior at the end of the rendezvous maneuver: the vertical error coordinate reaches the desired value without overshooting, thus ensuring the feasibility of the trajectory (i.e., it satisfies the constraints~\eqref{sys:DockingConstraints}), see Figure~\ref{fig:StraightLinePathcomparisonK_ez} and the zoom in Figure~\ref{fig:StraightLineZcomparisonK100} at about $95$sec. 

Finally, the sequence of rendezvous time is $82.9$sec, $61.8$sec, $52.4$sec, $46.52$sec for $k_{aggr} = 0.25, 0.5, 0.75, 1$, respectively. 
Comparing the rendezvous time with the desired one predicted by~\eqref{eq:t_des}, we observe a good matching expect for the case $k_{aggr}=1$. In fact, due to the transient behavior at the end of the aggressive maneuver, 
the optimal rendezvous time is $46.52$sec (note that $e_z = -0.1$ for $t = 96.52$sec, see Figure~\ref{fig:StraightLineZcomparisonK100}), yet the desired rendezvous time is $40.31$sec.

\subsection{Rendezvous with coupled longitudinal and lateral motion}
\begin{figure}[htpb]
  \begin{center}
     \includegraphics[width=9cm]{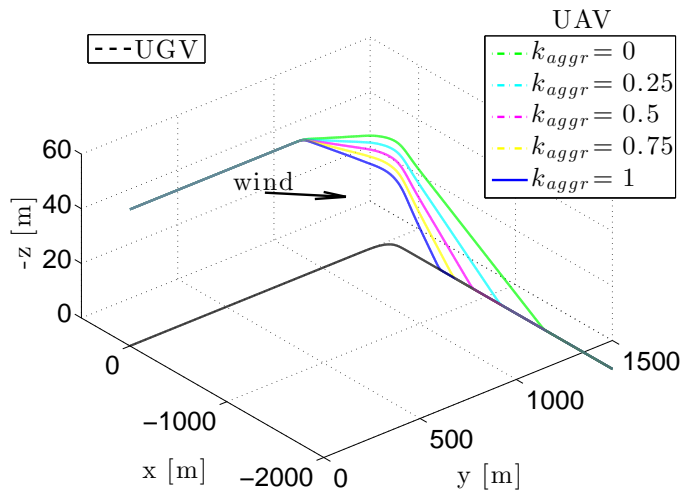}
\caption{Three dimension path of the UAV for a complex scenario for different $k_{aggr} = \{0, 0.25, 0.5, 0.75,1\}$.
    } \label{fig:90degPath}
  \end{center}
\end{figure}

\begin{center}
\begin{figure*}[!ht]
     \subfloat[${e}_z$.]
     {\includegraphics[width=5.5cm]{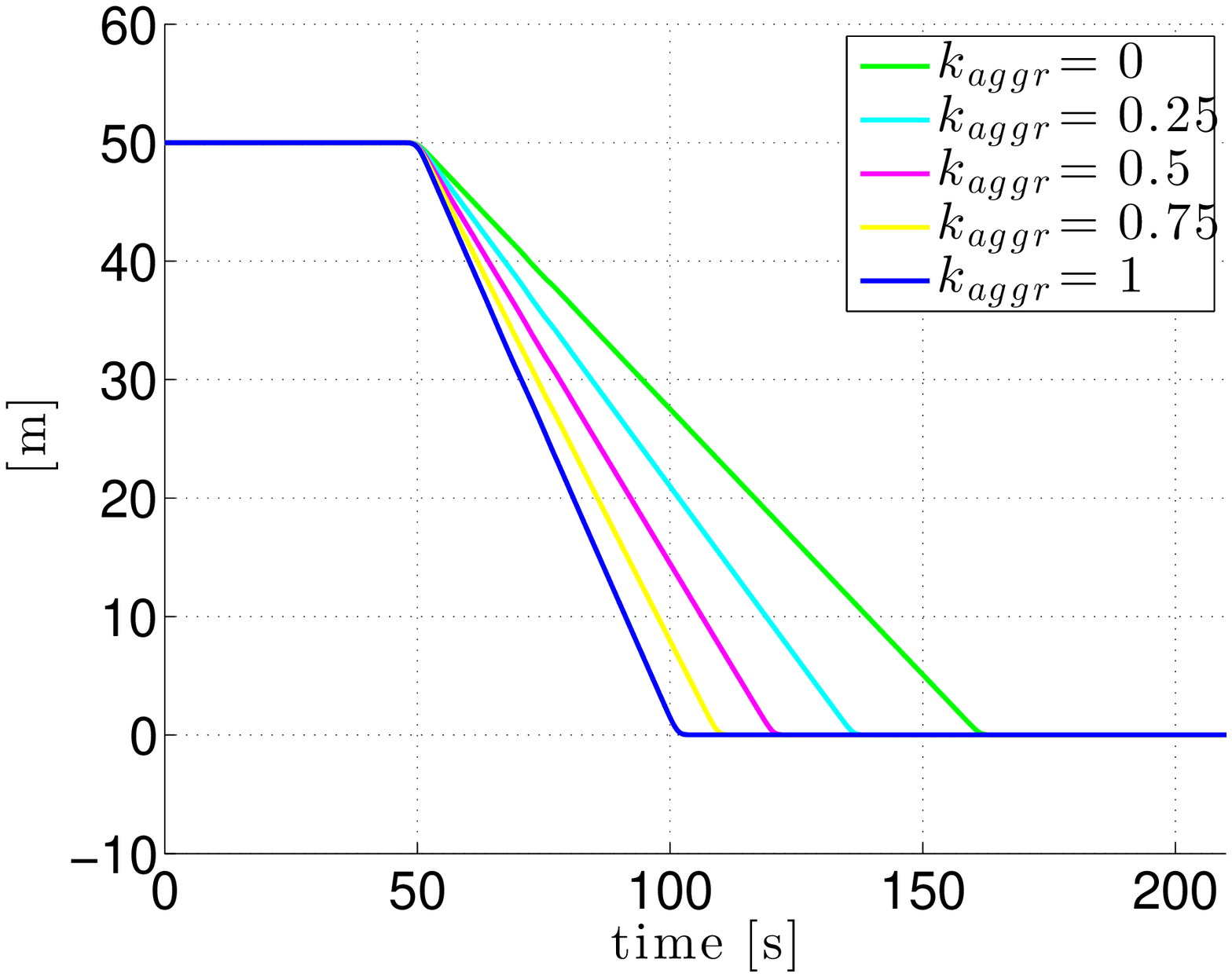} \label{fig:90deg_K_ez}}
     \subfloat[${e_v}$.]
     {\includegraphics[width=5.5cm]{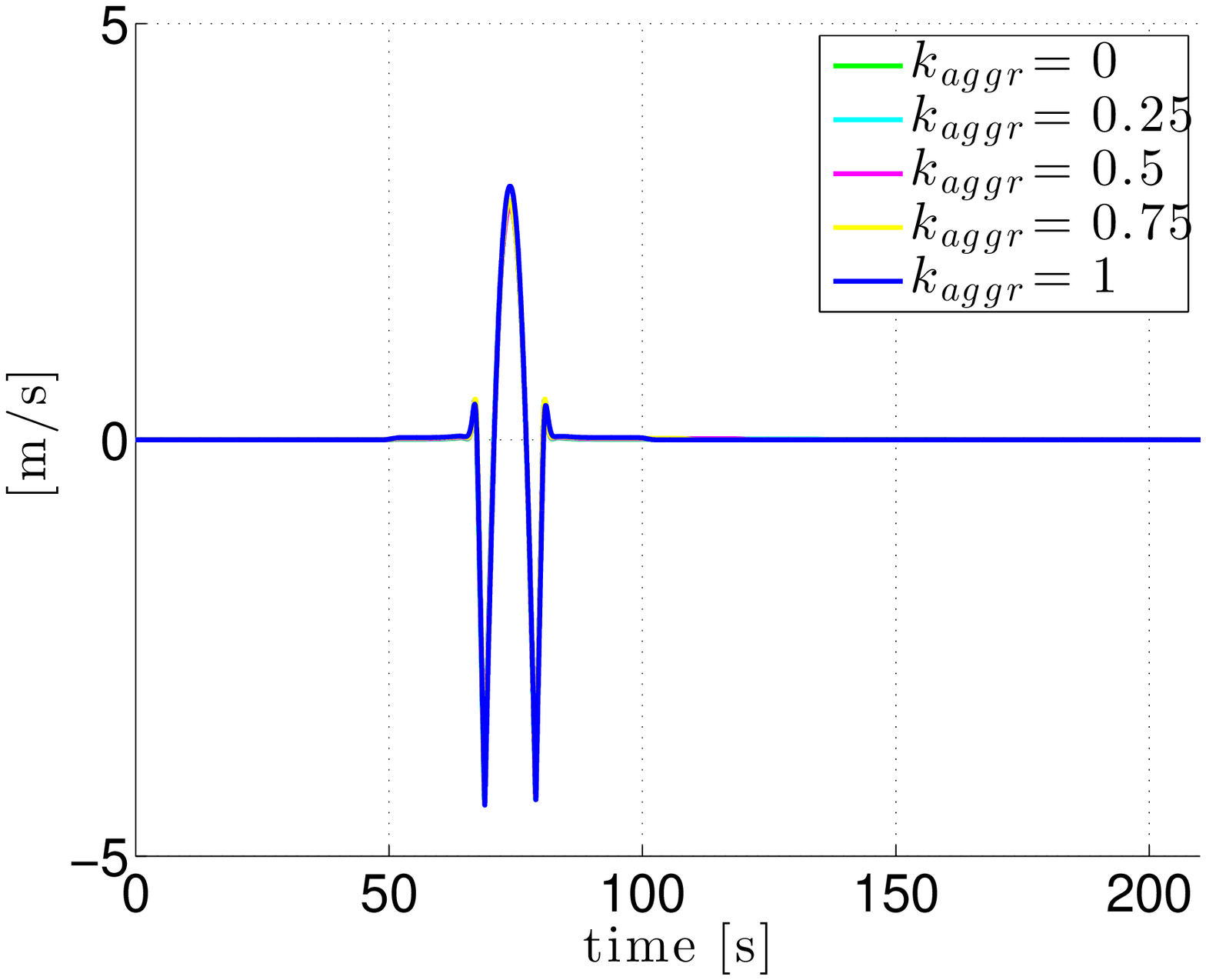}\label{fig:90deg_K_ev}}
     \subfloat[${e_\gamma}$.]
     {\includegraphics[width=5.5cm]{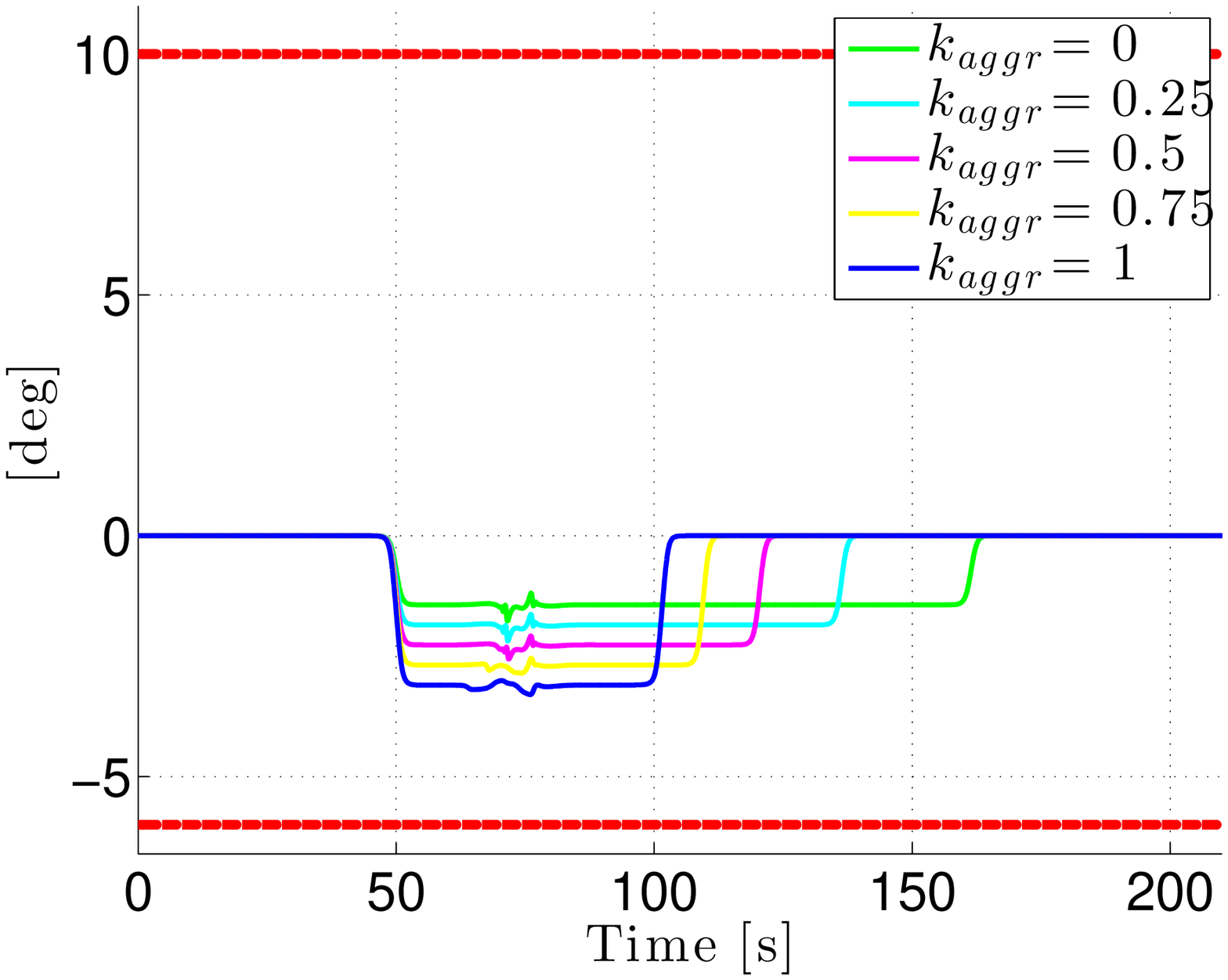} \label{fig:90deg_K_ez_gam}}

     \subfloat[${u_1}$.]
     {\includegraphics[width=5.5cm]{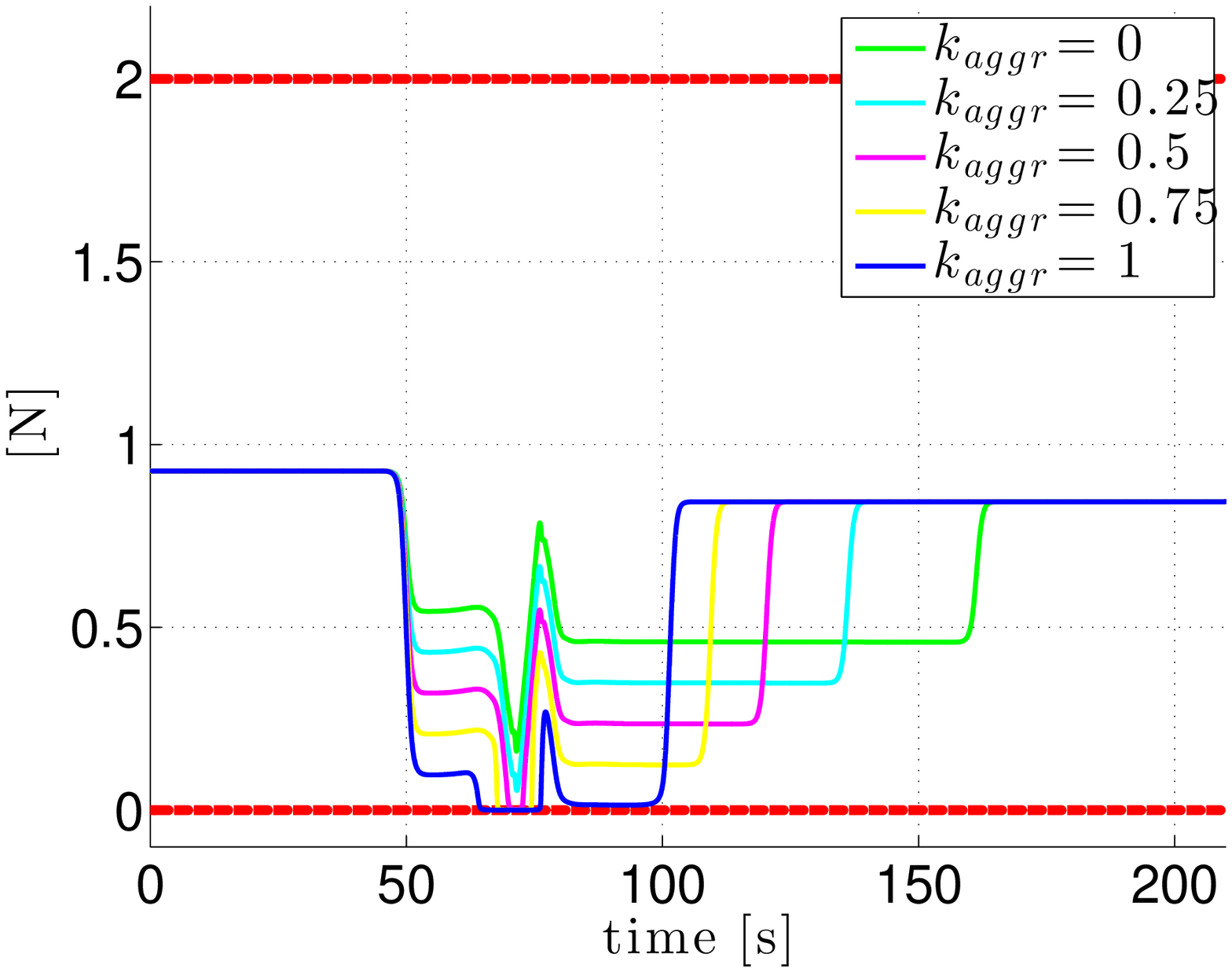} \label{fig:90deg_K_ez_u1}}
     \subfloat[${u_3}$.]
     {\includegraphics[width=5.5cm]{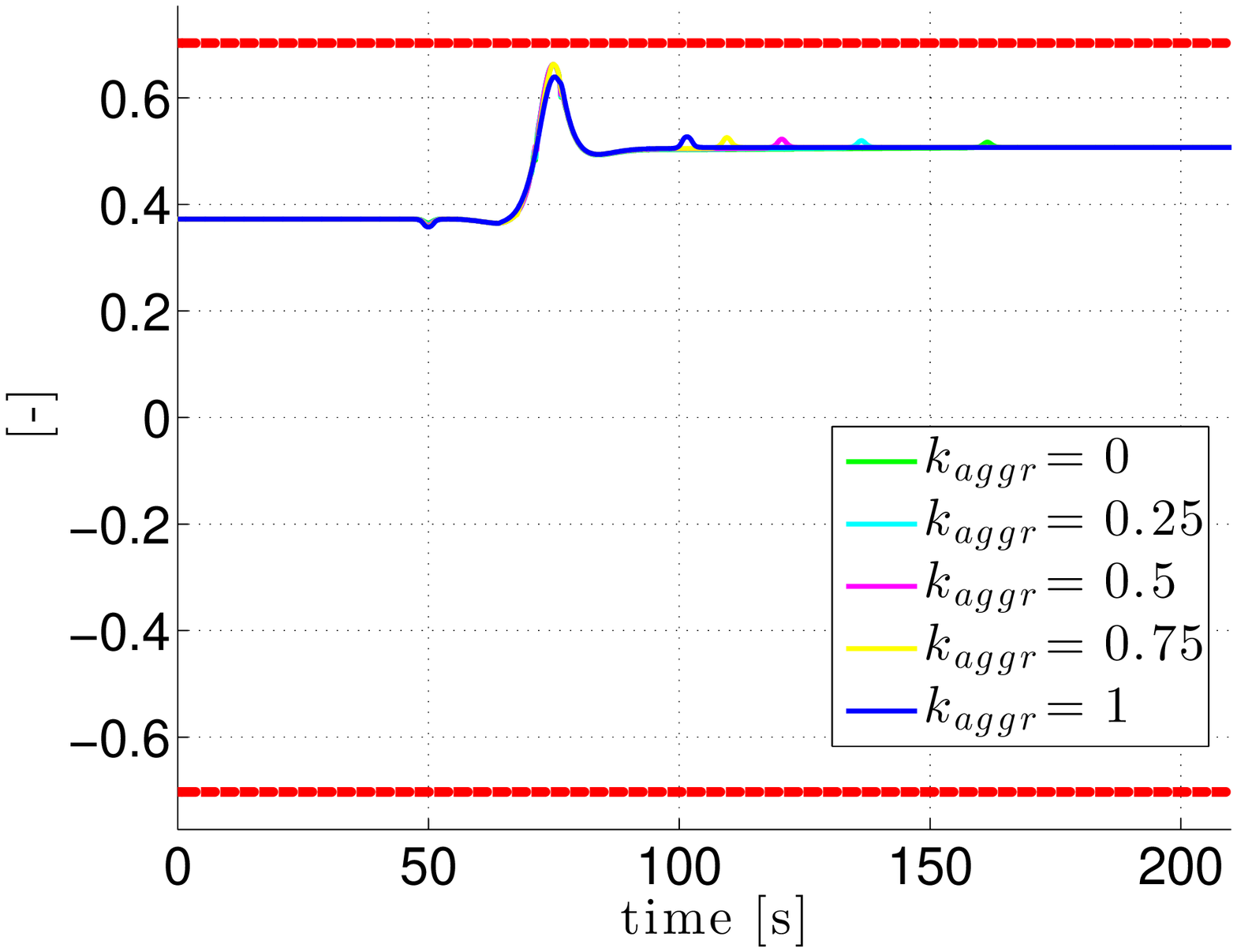} \label{fig:90deg_K_ez_CL}}
     \subfloat[$n_{lf}$.]
     {\includegraphics[width=5.5cm]{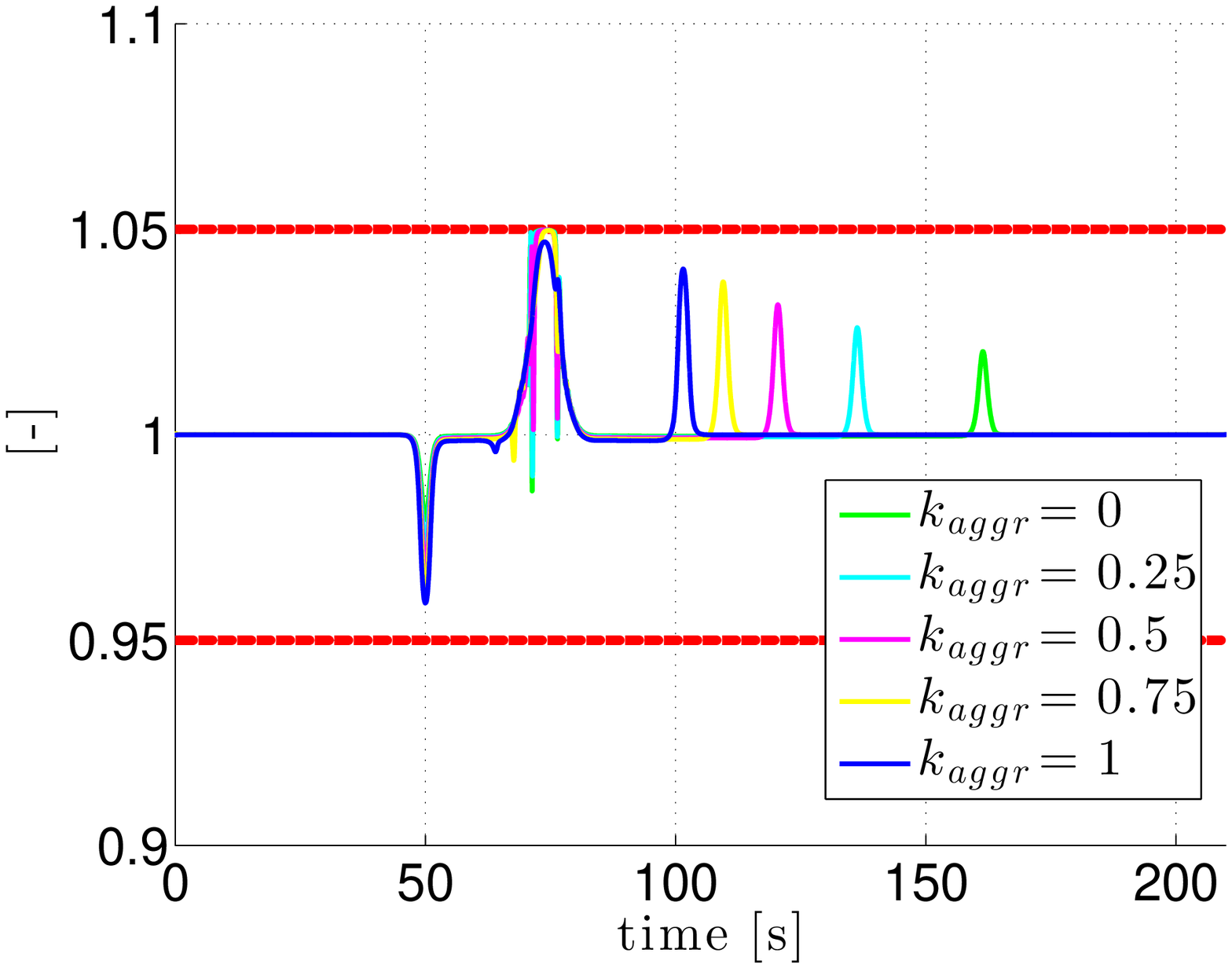} \label{fig:90deg_K_ez_nlf}}
     \caption{Rendezvous with coupled longitudinal and lateral motion (a) vertical error coordinate, 
     (b) error speed, 
     (c) error flight path angle, (d) thrust, (e) coefficient lift, and (f) load factor for $k_{aggr} = \{0, 0.25, 0.5, 0.75,1\}$. Constraints are in dashed line. 
     } \label{fig:90deg_K}
\end{figure*} 
\end{center}

In this scenario, the UGV is moving along a circuit as the one mentioned in the Introduction. 
In particular, we take into account a section of the circuit shown in Figure~\ref{fig:scenario2} which is composed by $90^\circ$ turn with a radius of $35$m and straights of $1200$m before and after the turn, see Figure~\ref{fig:90degPath}. 
The initial position of the UAV is $(x_A, y_A, z_A) = (0, 0, -50)$, the orientation is $\chi_A = 0$, and flight angle and roll angle are $\gamma_A = 0$, $\phi_A = 0$, respectively. 
The initial position and orientation of the UGV are $(x_G, y_G, z_G) = (0, 0, 0)$ and $\chi_G = 0$, respectively. 
We run the rendezvous trajectory generation strategy for aggressiveness index equals to $k_{aggr} = \{0, 0.25, 0.5, 0.75,1\}$ and we compare the (local) optimal trajectories in Figures~\ref{fig:90degPath} and~\ref{fig:90deg_K}. 
As in the previous computations, we are able to control the aggressiveness of the UAV trajectory. 
Thus, for $k_{aggr} = 1$ (blue line in Figures~\ref{fig:90degPath} and~\ref{fig:90deg_K}) several constraints are active. 
For $k_{aggr} = 0$, the UAV height is reduced gradually thus highlighting the soft feature of the local optimal trajectory (green line in Figures~\ref{fig:90degPath} and~\ref{fig:90deg_K}). 

\begin{center}
\begin{figure*}[!ht]
     \subfloat[${e_y}$.]
     {\includegraphics[width=5.5cm]{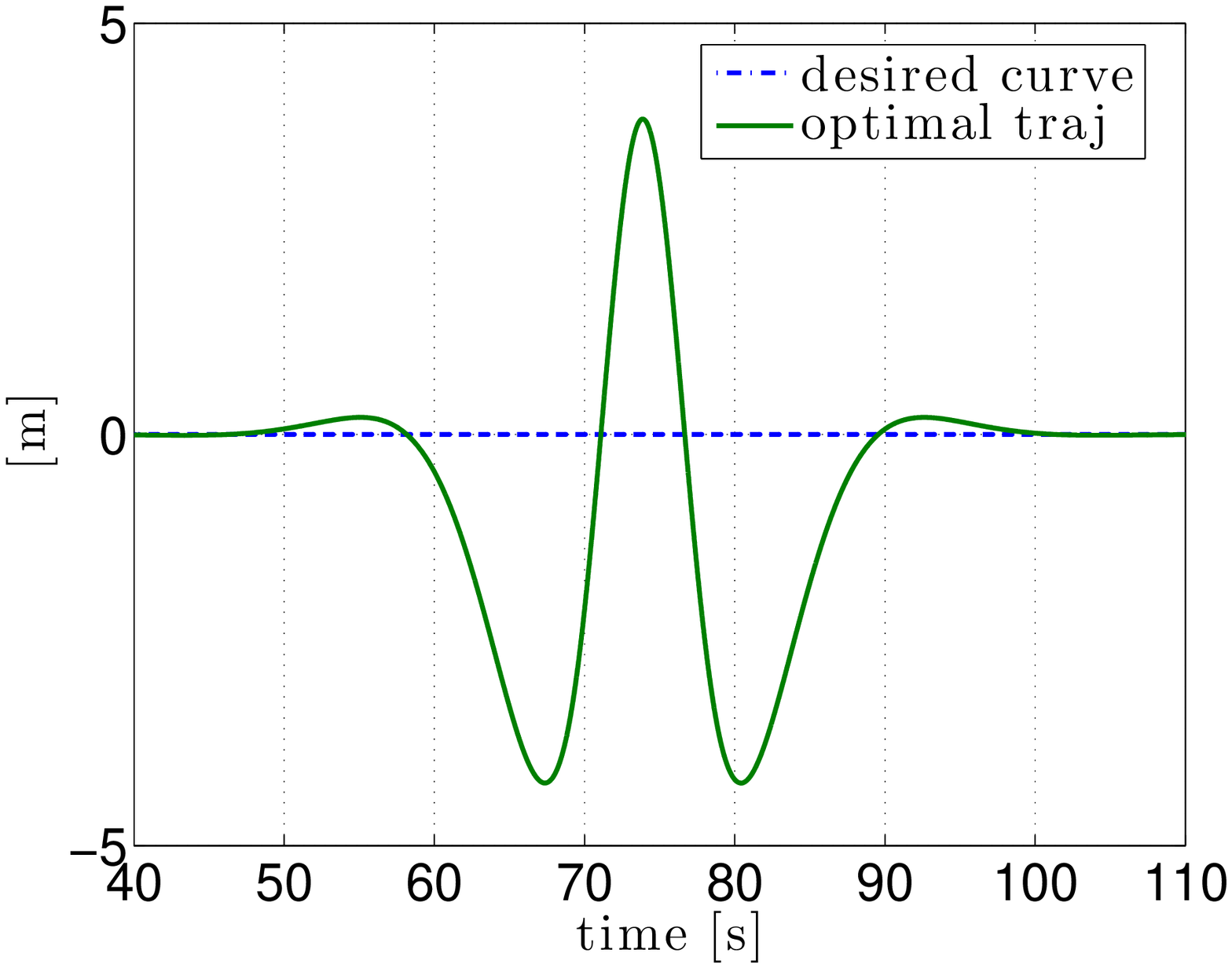}\label{fig:90degK1_ey}}
     \subfloat[${\phi_A}$.]
     {\includegraphics[width=5.5cm]{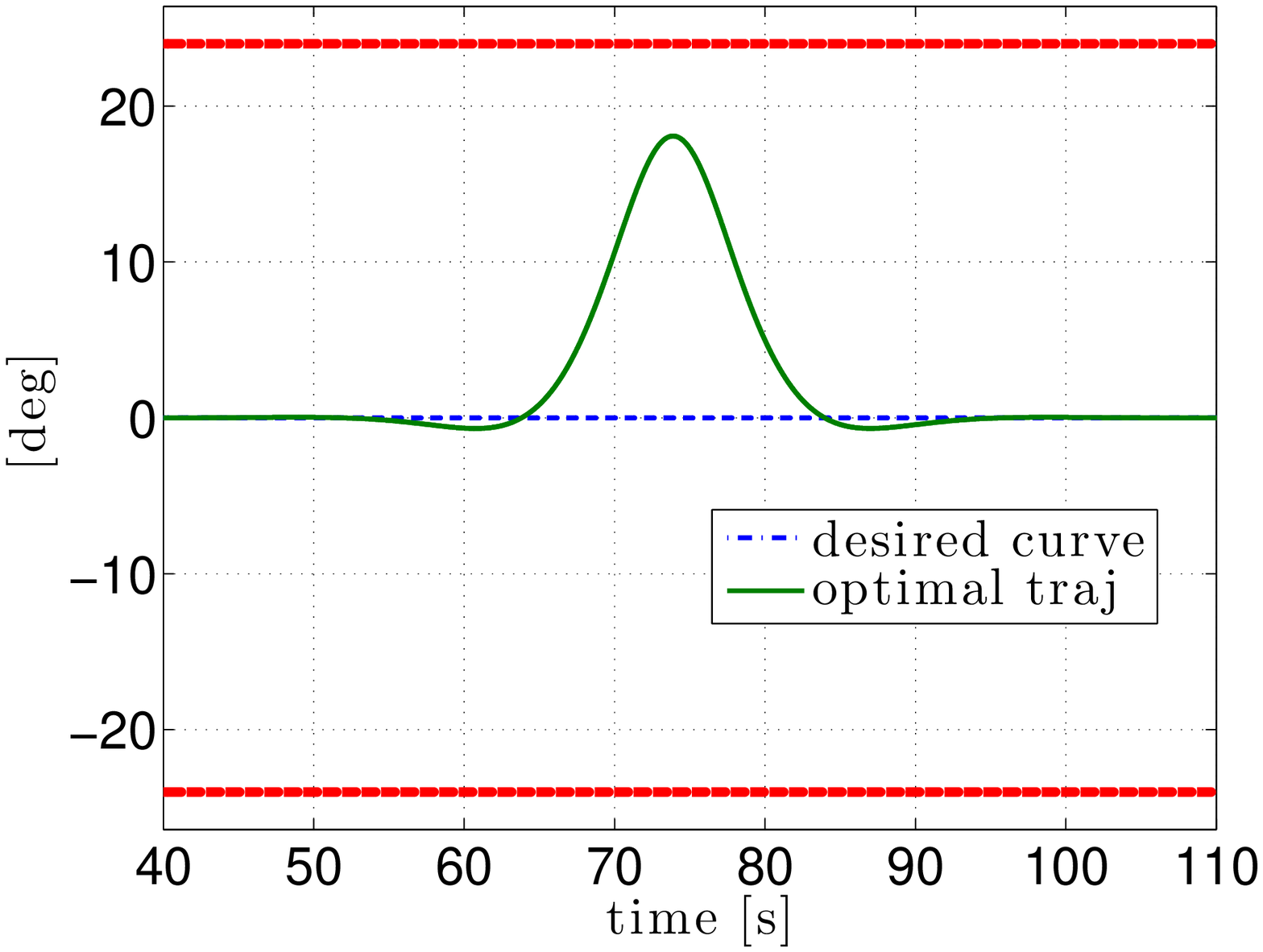} \label{fig:90degK1_ph}}
     \subfloat[$u_4$ vs $a_{lat}$.]
     {\includegraphics[width=5.5cm]{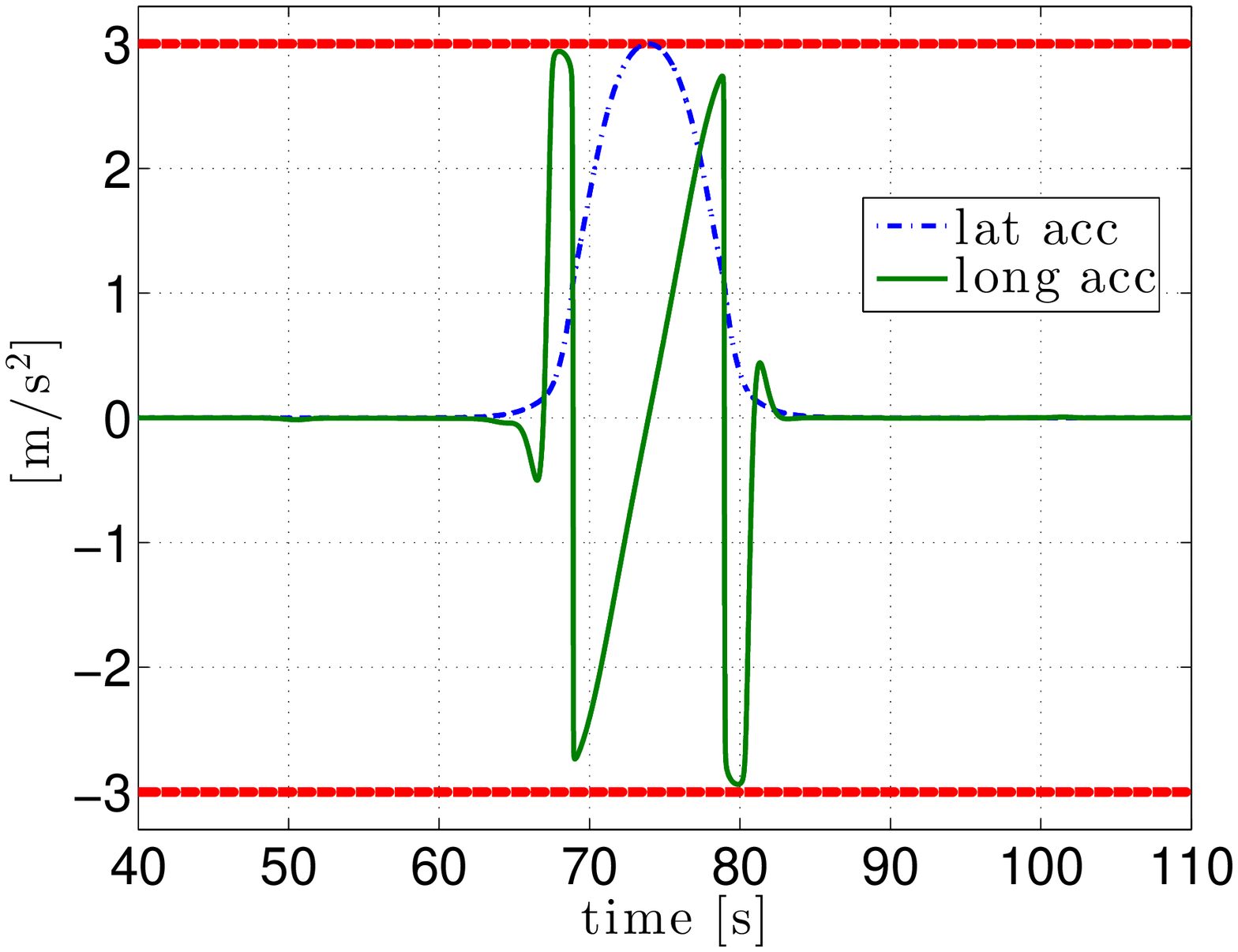} \label{fig:90degK1_acar}}     
     \caption{Rendezvous with coupled longitudinal and lateral motion for $k_{aggr} = 1$ (a) lateral error coordinate, (b) UAV roll angle, and (c) UGV longitudinal-lateral acceleration. Constraints are in dashed line. 
     } \label{fig:90degK1}
\end{figure*} 
\end{center}

It is worth highlighting the effect of the right turn on the rendezvous maneuver.  
In order to minimize the lateral error coordinate, the UAV turns by rolling, see Figure~\ref{fig:90degK1_ph}. 
However, the UAV is not able to track exactly the UGV. The lateral error coordinate is no zero and the constraint on the roll angle is never active, see Figures~\ref{fig:90degK1_ey} and~\ref{fig:90degK1_ph}. 
This is due to the fact that the constraint on the load factor becomes active (see Figure~\ref{fig:90deg_K_ez_nlf}) before the roll angle reaches its maximum value. Indeed, in constant descent flight conditions (i.e., $\dot{\gamma}_A=0$ and $\gamma_A < 0$), the lift must be equal to $m g \frac{\cos{\gamma_A}}{\cos{\phi_A}}$ and the load factor becomes $n_{lf} = \frac{\cos{\gamma_A}}{\cos{\phi_A}}$. 
It is evident that the roll angle is constrained by $\arccos{ \frac{\cos{\gamma_A}}{n_{lf \, max}} }$ which turns out to be less than $\phi_{max}$. This explains the no-zero lateral error coordinate. 

This computation allows us to highlight the coupled UAV-UGV dynamics. Indeed, for $k_{aggr} = 1$, the UAV roll reaches $18$deg (as discussed before, its maximum value) and the UGV lateral acceleration is at the maximum value, $a_{max}$, at exactly the same time $t=73.8$sec, see Figure~\ref{fig:90degK1_ph} and~\ref{fig:90degK1_acar}, respectively.  
Finally, although the desired curve is based on the decoupled UAV-UGV dynamics, we are able to predict the rendezvous time. 
The sequence of rendezvous time is $111.8$sec, $86.9$sec, $71.1$sec, $60.3$sec, $52.4$sec for $k_{aggr} = 0, 0.25, 0.5, 0.75, 1$, respectively, as predicted by~\eqref{eq:Trendezvous}.

\section{Conclusions} \label{sec:conclude}
In this paper we proposed an optimal control approach for the refueling problem of fixed-wing UAVs using a UGV as a refueling unit. 
We provided a rigourous optimal control problem formulation for UAV rendezvous with the moving UGV and 
we addressed the optimal control problem by using a trajectory-tracking approach. Based on a nonlinear optimal control solver, we proposed an optimal control based strategy which allows us to compute optimal feasible trajectories for both UAV and UGV. 
By changing the aggressiveness index in our proposed strategy, we are able to compute aggressive trajectories (i.e., several constraints are active while the UAV is approaching the UGV) or very smooth ones. 
A key property of the proposed approach is that we are able to predict and, therefore, select (in form of tuning knob) the time to rendezvous, which is an important performance feature of the UAV trajectory. 
We provided numerical computations showing the effectiveness of the proposed approach. 
Future directions of research will include field tests where the obtained optimal trajectories are feed as reference trajectories to the trajectory tracking algorithms that are running on the vehicles to perform the rendezvous task. 

\appendix 

\label{sec:appendix}
The UAV parameters are based on the ``Zagi'' flying wing~\cite{beard2012small}: 
\[
  m = 1.56\, \text{kg}
  \quad
  S = 0.2589 \, \text{m$^2$} 
  \quad
  b = 1.4224 \, \text{m} 
\]
\[
  C_{D0} = \, 0.01631  
  \quad  
  k_{D/L} = \, 0.04525  
\]
We assume that the air density $\rho$ is constant and equal to $1.225$[kg/m$^3$]. 
The minimum and maximum airspeed, normal load, maximum thrust, roll angle, and coefficient lift are set as follows
\[
  v_{min} = 12\,\text{m/s} \,,
  \quad
  v_{max} = 20\,\text{m/s}\,,
\]
\[
  n_{lf\, min} = 0.95 \,,
  \quad
  n_{lf\, max} = 1.05 \,,
\]
\[
  \gamma_{min} = -6\,\text{deg}\,,
  \quad
  \gamma_{max} = 10 \,\text{deg}\,,
\]
\[
  \phi_{max} = 24\, \text{deg}\,,
  \quad
  u_{1 max} = 2\, \text{N}\,,
\]
\[
  u_{2 max} = 5 \,\text{deg/s}\,,
  \quad
  u_{3 max} = 0.7\,.
\]
The maximum acceleration of the UGV is 
\[
  a_{max} = 3 \, \text{m/s$^2$}\,.
\]
The maximum course angle is defined by 
\[
  \bar{e}_{x} = \bar{e}_{y} = \bar{e}_{z} = 30\,\text{m}\,,
  \quad
  \bar{e}_{\chi} = 2\,\text{deg}\,.
\]

\ifCLASSOPTIONcaptionsoff
  \newpage
\fi



\end{document}